\documentclass[11pt]{amsart2000}
\usepackage{amsmath2000,amsthm2000,amsfonts,amscd2000,amssymb,eucal,latexsym}

\theoremstyle{definition}
\newtheorem{theorem}{Theorem}[section]
\newtheorem{corollary}[theorem]{Corollary}
\newtheorem{lemma}[theorem]{Lemma}
\newtheorem{proposition}[theorem]{Proposition}

\newtheorem{notation}[theorem]{Notation}
\newtheorem{example}[theorem]{Example}
\newtheorem{question}[theorem]{Question}
\newtheorem{definition}[theorem]{Definition}

\newcommand{\diam}{\text{\rm diam}}
\newcommand{\dist}{\text{\rm dist}}

\newcommand{\id}{\text{\rm id}}
\newcommand{\Aut}{\text{\rm Aut}}
\newcommand{\re}{\text{\rm Re}}
\newcommand{\im}{\text{\rm Im}}

\newcommand{\sa}{\text{\rm sa}}
\newcommand{\Rcp}{\text{\rm Rcp}}
\newcommand{\CPA}{CPA}

\newcommand{\rank}{\text{\rm rank}}

\newcommand{\Afn}{\text{\rm Afn}}

\newcommand{\op}{\text{\rm op}}
\newcommand{\fa}{\text{\rm fa}}
\newcommand{\alg}{\text{\rm alg}}

\allowdisplaybreaks

\begin{document}

\title{Matricial quantum Gromov-Hausdorff distance}

\author{David Kerr}
\address{Graduate School of Mathematical Sciences,
University of Tokyo,
3-8-1 Komaba, Meguro-ku, Tokyo 153-8914, Japan}
\email{dkerr@ms.u-tokyo.ac.jp}
\date{July 28, 2002}

\begin{abstract}
We develop a matricial version of Rieffel's Gromov-Hausdorff
distance for compact quantum metric spaces within the setting 
of operator systems and unital $C^*$-algebras. Our approach yields
a metric space of ``isometric'' unital complete
order isomorphism classes of metrized operator systems which in many cases
exhibits the same convergence properties as those in the quantum metric 
setting, as for example in Rieffel's approximation of the sphere by matrix 
algebras using Berezin quantization. Within the metric subspace of metrized 
unital $C^*$-algebras we establish the convergence of sequences which
are Cauchy with respect to a larger Leibniz distance, and we also prove an 
analogue of the precompactness theorems of Gromov and Rieffel.
\end{abstract}

\maketitle

\section{Introduction}

A compact quantum metric space, as defined by Rieffel in \cite{GHDQMS},
is an order-unit space equipped with a certain type of semi-norm,
called a Lip-norm, which plays the role of a Lipschitz semi-norm on
functions over a compact metric space. The crucial part of the definition of 
a Lip-norm $L$ on an order-unit space $X$ is the requirement that the metric
$$ \rho_L (\mu , \nu ) = \sup \{ | \mu (a) - \nu (a) | : 
a\in A \text{ and }L(a)\leq 1 \} $$
on the state space of $X$ give rise to the weak$^*$ topology. By applying 
Hausdorff distance to state spaces, Rieffel defines a quantum analogue of 
Gromov-Hausdorff distance and thereby synthetically obtains a complete 
separable metric space of ``isometric'' order isomorphism classes of compact 
quantum metric spaces for which a Gromov-type precompactness theorem 
holds \cite{GHDQMS}. The most immediate motivation for 
introducing a theory of quantum Gromov-Hausdorff distance is the search for 
an analytic framework for describing, or at least clarifying at a 
metric level, the type of convergence of spaces that has recently begun
to play a central role in string theory 
(see the introduction to \cite{GHDQMS} for a discussion and
references). The main objects of study thus tend to be $C^*$-algebras,
and so it is natural to ask, as does Rieffel in \cite{GHDQMS}, if it is
possible to develop a matricial version of quantum Gromov-Hausdorff
distance. This is the goal of the present paper.

The key is to define metrics on matrix states spaces using a Lip-norm
just as one does for an order-unit state space as above, only now replacing 
the modulus by matrix norms. We introduce this definition within a general
operator system setting in Section~\ref{S-Lip-norms}. We then define 
``complete'' distance (Section~\ref{S-distance}) by using Hausdorff distance 
at the matrix state space level in the same way that Rieffel does with regard 
to order-unit state spaces in the formulation of quantum Gromov-Hausdorff 
distance. In fact many of the constructions and arguments involving 
quantum Gromov-Hausdorff distance in \cite{GHDQMS,MACS} are naturally 
suited to our matricial setting and lead to similar estimates, as for
instance in the proof of the triangle inequality 
(Proposition~\ref{P-triangle}) and the approximation of the sphere
by matrix algebras via Berezin quantization (Example~\ref{E-sphere}). 
On the other hand a completely 
different approach is required to show that complete
distance zero implies ``isometric'' unital complete order isomorphism
(the subject of Section~\ref{S-zero}), and the proofs of the 
convergence of sequences of metrized unital $C^*$-algebras which 
are Cauchy with respect to ``$f$-Leibniz complete'' distance
(Section \ref{S-complete}) and our analogue of the Gromov and Rieffel
precompactness theorems (Sections \ref{S-tb}) ultimately 
rely on some arguments especially attuned to the complete order
context. 

This work was supported by the Natural 
Sciences and Engineering Research Council of Canada. I thank 
Yasuyuki Kawahigashi and the operator algebra group at the
University of Tokyo for their hospitality and for the
invigorating research environment they have provided. I also thank
Hanfeng Li for making a careful critical reading of an initial draft that 
prompted a number of clarifications and corrections.

\section{Lip-normed operator systems and matrix
state space metrization}\label{S-Lip-norms}

We begin by describing our operator system 
framework. For references see \cite{ER,Pau,Was}. 
A (concrete) {\em operator system} is a closed unital self-adjoint linear
subspace of a unital $C^*$-algebra (for an abstract definition see
\cite{ER}). Given an operator 
system $X$, for each $r > 0$ we will denote by $\mathcal{B}^X_r$ the 
closed norm ball $\{ x\in X : \| x \| \leq r \}$ of radius $r$. The state 
space of $X$ will be denoted by $S(X)$. We will denote by $X_{\sa}$ the set 
of self-adjoint elements of $X$. The unit of $X$ will be written
$1$, or sometimes $1_X$ for clarity. For $x\in X$ we write $\re (x)$ 
and $\im (x)$ to refer to the self-adjoint elements 
$(x+x^* )/2$ and $(x-x^* )/(2i)$ (the real and 
imaginary parts of $x$), respectively.

The $C^*$-algebra of $n\times n$ matrices over $\mathbb{C}$ will be
written $M_n$. Given operator systems $X$ and $Y$ 
we say that a linear map $\varphi : X \to Y$ is {\em $n$-positive}
if the map $\id_n \otimes\varphi : M_n \otimes X \to M_n \otimes Y$
is positive, and if $\id_n \otimes\varphi$ is $n$-positive 
for all $n\in\mathbb{N}$ then we say
that $\varphi$ is {\em completely positive}. A completely positive (resp.\
unital completely positive) linear map will be referred to as a 
{\em c.p.}\ (resp.\ {\em u.c.p.}) map. If $\varphi : X \to Y$ is a unital 
$m$-positive map with $m$-positive inverse for $m=1, \dots , n$ then
$\varphi$ is a {\em unital $n$-order isomorphism}, and if $\varphi$ is 
u.c.p.\ with c.p.\ inverse then $\varphi$ is a {\em unital complete
order isomorphism}.

An operator system $X$ is {\em nuclear} if the identity
map on $X$ lies in the point-norm closure of the set of u.c.p.\ maps
from $X$ to itself which factor through matrix algebras.

Given an operator system $X$ and $n\in\mathbb{N}$, there is a bijective
linear map from c.p.\ maps $X\to M_n$ to positive linear functionals
on $M_n \otimes X$ \cite[Thm.\ 5.1]{Pau} defined as follows. To each
c.p.\ map $\varphi : X\to M_n$ we associate the positive linear functional
$\sigma_\varphi$ on $M_n \otimes X$ given by
$$ \sigma_\varphi ((x_{ij})) = \frac1n \sum_{1\leq i,j \leq n} \varphi
(x_{ij})_{ij} $$
for all $(x_{ij})\in M_n (X) \cong M_n \otimes X$.
Conversely, to each positive linear functional $\sigma$ on 
$M_n \otimes X$ we associate the c.p.\ map $\varphi_\sigma : X \to M_n$
given by
$$ (\varphi_\sigma (x))_{ij} = n\,\sigma (e_{ij} \otimes x) $$
for all $x\in X$, 
where $\{ e_{ij} : 1\leq i,j \leq n \}$ is the set of standard matrix 
units of $M_n$. The maps $\varphi\mapsto\sigma_\varphi$ and 
$\sigma\mapsto\varphi_\sigma$ are mutual inverses and are homeomorphisms with 
respect to the point-norm topologies (for the space of positive linear
functionals this is the weak$^*$ topology) as well as with respect
to the norm topologies. If $\varphi : X\to M_n$ is u.c.p.\ then
$\sigma_\varphi$ is a state on $M_n \otimes X$. However, if $\sigma\in
S(M_n \otimes X)$ then $\varphi_\sigma$ need not be unital, nor even
contractive, although it is clear that $\| \varphi_\sigma \| \leq n^3$ 
(see the discussion after Theorem 5.4 in \cite{Pau}). We denote by 
$SCP_n (X)$ the collection of c.p.\ maps $\varphi : X \to M_n$ such that
$\sigma_\varphi$ is a state on $M_n \otimes X$, and by 
$UCP_n (X)$ the subcollection of $SCP_n (X)$ consisting of all
u.c.p.\ maps from $X$ into $M_n$ (the {\em matrix state spaces}). 

We now introduce metrics into our picture via 
the notion of a Lip-norm, which we recall from \cite{GHDQMS}.

\begin{definition}[{\cite[Defns. 2.1 and 2.2]{GHDQMS}}]\label{D-cqms}
Let $A$ be an order-unit space. A {\em Lip-norm} on $A$ 
is a semi-norm $L$ on $A$ such that
\begin{enumerate}
\item[(1)] for all $a\in A$ we have $L(a) = 0$ if and only if $a$ is
a scalar multiple of the order unit, and
\item[(2)] the metric $\rho_L$ defined on the state space $S(A)$ by
$$ \rho_L (\mu , \nu ) = \sup \{ |\mu (a) - \nu (a) | : 
a\in A \text{ and }L(a)\leq 1 \} $$
induces the weak$^*$ topology.
\end{enumerate}
A pair $(A,L)$ consisting of an order-unit space $A$ with Lip-norm $L$
is called a {\em compact quantum metric space}.
\end{definition}

Important examples of order-unit spaces are real linear unital subspaces of 
self-adjoint elements in an operator system, and in fact every order-unit
space is isomorphic to one of these, as shown in Appendix 2 of 
\cite{GHDQMS}. We can thus apply the above definition in a direct way
to our setting. First we introduce some general notation.

\begin{notation}\label{N-L}
Let $X$ be an operator system and $L$ a semi-norm on a linear subspace of 
$X$ or a real linear subspace of $X_{\sa}$. We denote by $\mathcal{D}(L)$
the domain of $L$, or, if $L$ is permitted to take the value $+\infty$,
the set of elements in the domain of $L$ on which $L$ is finite-valued. 
For $r>0$ we denote by
$\mathcal{D}_r (L)$ the set $\{ x\in\mathcal{D}(L) : L(x) \leq r \}$.
\end{notation}

\begin{definition}\label{D-clipnorm}
By a {\em Lip-normed operator system} we mean a pair $(X,L)$ where
$X$ is an operator system and $L$ is a Lip-norm on a dense order-unit
subspace of $X_{\sa}$ such that $\mathcal{D}_1 (L)$ is closed in 
$X_{\sa}$. If $X$ is a unital $C^*$-algebra then we will
also refer to $(X,L)$ as a {\em Lip-normed unital $C^*$-algebra}. Any
qualifiers preceding ``Lip-normed'' will refer to the Lip-norm while
those following it will refer to the operator system or $C^*$-algebra
(e.g., {\em lower semicontinuous Lip-normed nuclear operator system}). 
\end{definition}

A Lip-norm $L$ on an order-unit space $A$ is said to be {\em closed} if 
if the set $\{ a\in A : L(a) \leq 1 \}$ is closed in the completion of $A$
\cite[Defn. 4.5]{MSS}. Thus the requirement in Definition~\ref{D-clipnorm} 
that $\mathcal{D}_1 (L)$ be closed in $X_{\sa}$ is equivalent to asking that 
$L$ be a closed Lip-norm. Given any Lip-norm 
$L$ on an order-unit space $A$ there is
a largest lower semicontinuous Lip-norm $L^s$ smaller than $L$
\cite[Thm.\ 4.2]{MSS}, and $L^s$ extends to a closed Lip-norm $L^c$
\cite[Prop.\ 4.4]{MSS}.
The theorem and proposition from \cite{MSS}
cited in the last sentence also show that $\rho_{L^c} = \rho_{L^s} = \rho_L$. 
Furthermore, the property of being closed passes to order-unit quotients by
\cite[Prop.\ 3.3]{GHDQMS}, and it also holds in natural examples 
of interest---see for instance Example~\ref{E-actions} 
and \cite[Prop.\ 3.6]{MSS}. Thus, in view of the completeness of
operator systems, it is natural to assume that our Lip-norms are
closed. This will also guarantee that complete distance zero is
equivalent to the existence of a unital complete order isomorphism
which is {\em bi-Lip-isometric} in the obvious sense:

\begin{definition}
Let $(X,L_X )$ and $(Y,L_Y )$ be Lip-normed operator systems. We will say
that a positive unital map $\Phi : X \to Y$ is {\em Lip-isometric} if 
$\Phi (\mathcal{D}(L_X ))\subset\mathcal{D}(L_Y )$ and $L_Y (\Phi (x)) 
= L_X (x)$ for all $x\in\mathcal{D}(L_X )$. If $\Phi$ has a positive 
inverse then we say that $\Phi$ is {\em bi-Lip-isometric} if both $\Phi$ and 
$\Phi^{-1}$ are Lip-isometric.
\end{definition}

If we were to define a strict operator system analogue of a Lip-norm
then the conditions on the semi-norm $L$ in the following proposition 
would seem to be the most reasonable. Indeed many examples
arise naturally in this way, as Example~\ref{E-actions} illustrates.

\begin{proposition}\label{P-selfadj}
Let $L$ be a semi-norm on an operator system $X$, permitted to take the
value $+\infty$, such that $\mathcal{D}(L)$ is dense in $X$, 
$\mathcal{D}_1 (L)$ is closed in $X_{\sa}$, and
\begin{enumerate}
\item[(i)] $L(x^* ) = L(x)$ for all $x\in A$ (adjoint invariance),
\item[(ii)] for all $x\in X$ we have $L(x) = 0$ if and only if 
$x\in\mathbb{C}1$, and
\item[(iii)] the metric $d_L (\sigma , \omega ) = 
\sup_{x\in\mathcal{D}_1 (L)} | \sigma (x) - \omega (x) |$ on $S(X)$ induces 
the weak$^*$ topology.
\end{enumerate}
Then the restriction $L'$ of $L$ to the order-unit space 
$\mathcal{D}(L) \cap X_{\sa}$ is a Lip-norm, $(X,L' )$ is a 
Lip-normed operator system, and the restriction map from $S(X)$ onto
$S(\mathcal{D}(L' ))$ is a weak$^*$ homeomorphism which is isometric for 
$d_L$ and $\rho_{L'}$.
\end{proposition}

\begin{proof}
First note that the fact that $\mathcal{D}(L)$ is closed in $X$
immediately implies that $\mathcal{D}(L' )$ is closed in $X_{\sa}$. Next, 
if $x\in X_{\sa}$ and $\epsilon > 0$ then by the density of $\mathcal{D}(L)$
we can find a $y\in\mathcal{D}(L)$ with $\| x - y \| < \epsilon$. Then
$$ \| x - \re (y) \| \leq \| x-y \| /2 + \| (x-y)^* \| /2 <
\epsilon $$
while $L' (\re (y)) = (L(y) + L(y^* ))/2 = L(y) < +\infty$, and so
$\mathcal{D}(L' )$ is dense in $X_{\sa}$. With this fact it is 
straightforward to show that the restriction map from $S(X)$ onto 
$S(\mathcal{D}(L' ))$ is a weak$^*$ homeomorphism. To see that this map
is isometric, suppose $\sigma , \omega\in S(X)$ and $\epsilon > 0$. Then
we can find an $x\in\mathcal{D}_1 (L)$ such that $| \sigma (x) - \omega 
(x) | \geq d_L (\sigma , \omega ) - \epsilon$, and so for some
complex number $\mu$ of unit modulus we have $\sigma (\mu x) - \omega
(\mu x) \geq d_L (\sigma , \omega ) -\epsilon$. Since
$$ L' (\re (\mu x)) \leq (L(\mu x ) + L(\bar{\mu}x^* ))/2 = 
(L(x) + L(x^* )) / 2 = L(x)\leq 1  $$
we then have
\begin{align*}
\rho_{L'}(\sigma |_{\mathcal{D}(L' )} , \omega |_{\mathcal{D}(L' )}) &\geq
\sigma (\re (\mu x)) - \omega (\re (\mu x)) \\
&\geq d_L (\sigma , \omega ) -\epsilon ,
\end{align*}
from which we infer that the map in question is indeed an isometric 
weak$^*$ homeomorphism. Since condition (ii) in the proposition statement 
immediately implies condition (1) in Definition~\ref{D-cqms} for $L'$, it 
thus follows that $L'$ is a Lip-norm, and so $(X,L' )$ is a 
Lip-normed operator system in view of the density of 
$\mathcal{D}(L' )$ in $X_{\sa}$.
\end{proof}

In the converse direction, given a Lip-normed operator system $(X,L' )$ we 
can extend $L'$ to a semi-norm $L$ on $X$ such that the conditions and 
conclusions in the statement of Proposition~\ref{P-selfadj} hold. 
Definition~\ref{D-semi-norm} and Proposition~\ref{P-restriction} indicate 
how this can be done.

\begin{example}[ergodic actions of compact groups]\label{E-actions}
As studied in \cite{MSACG}, ergodic actions of compact groups give rise 
to important examples of Lip-normed unital $C^*$-algebras, notably  
noncommutative tori (see Example~\ref{E-tori}). 
Let $\gamma$ be an ergodic action of a compact group 
$G$ on a unital $C^*$-algebra $A$. Let $e$ be the identity element 
of $G$. We suppose that $G$ is 
equipped with a length function $\ell$, that is, a continuous function 
$\ell : G \to \mathbb{R}_{\geq 0}$ such that, for all $g,h\in G$,
\begin{enumerate}
\item[(1)] $\ell (gh)\leq \ell (g) + \ell (h)$, 
\item[(2)] $\ell (g^{-1}) = \ell (g)$, and 
\item[(3)] $\ell (g) = 0$ if and only if $g\neq e$. 
\end{enumerate}
The group action $\gamma$ and the length function $\ell$ 
together yield the semi-norm $L$ on $A$ defined by
$$ L(a) = \sup \left\{ \| \gamma_g (a) - a \| / \ell (g) : 
g\neq e \right\} . $$
It is easily seen that $L$ is adjoint-invariant and that
$L(a) = 0$ if and only if $a\in\mathbb{C}1$.
Furthermore, by \cite[Thm. 2.3]{MSACG}
the metric $d_L (\sigma , \omega ) = \sup_{x\in\mathcal{D}_1 (L)} 
| \sigma (x) - \omega (x) |$ on $S(A)$ induces 
the weak$^*$ topology, by \cite[Prop.\ 2.2]{MSACG} $\mathcal{D}(L)$ is 
dense in $X$, and it is readily verified that $\mathcal{D}_1 (L)$ is closed
in $A$ (see \cite[Prop.\ 8.1]{GHDQMS}), so that by 
Proposition~\ref{P-selfadj} we obtain a Lip-normed unital $C^*$-algebra by 
restricting $L$ to $\mathcal{D}(L) \cap A_{\sa}$.
\end{example}

\begin{example}[quotients]\label{E-quotients}
Let $(X,L)$ be a Lip-normed operator system, $Y$ an operator system,
and $\Phi : X \to Y$ a surjective unital positive linear map. 
Then by \cite[Prop.\ 3.1]{GHDQMS} $L$ gives rise to a Lip-norm $L_Y$ on 
$\Phi (\mathcal{D}(L))$ via the prescription
$$ L_Y (y) = \inf \{ L(x) : x\in \mathcal{D}(L)\text{ and } \Phi (x) = y 
\} $$
for each $y\in Y$, and the induced map from $(S(Y), \rho_{L_Y})$
to $(S(X), \rho_L )$ is an isometry. Since $\Phi (\mathcal{D}(L))$ is 
clearly dense in $Y_{\sa}$ and $L_Y$ is closed by \cite[Prop.\ 3.3]{GHDQMS}, 
$(Y,L_Y )$ is a Lip-normed operator system. We say that {\em $L$ induces 
$L_Y$ via $\Phi$}.
\end{example}

The following definition captures the observation that Lip-norms
define metrics on matrix state spaces in much the same way as they do
on state spaces. We will thereby be able to define a matrix
version of quantum Gromov-Hausdorff distance by applying 
Hausdorff distance to the matrix state spaces (Definition~\ref{D-dist}).

\begin{definition}\label{D-ucpmetric}
Let $(X,L)$ be a Lip-normed operator system and $n\in\mathbb{N}$. 
We define the metric $\rho_{L,n}$ on $UCP_n (X)$ by 
$$ \rho_{L,n} (\varphi , \psi ) = \sup_{x\in\mathcal{D}_1 (L)}
\| \varphi (x) - \psi (x) \| $$
for all $\varphi , \psi\in UCP_n (X)$,
\end{definition}

Note that $\rho_{L,n}$ is indeed a metric since it clearly satisfies the
triangle inequality and is symmetric, and it is non-zero at any pair
of distinct points owing to the density of $\mathcal{D}(L)$ in $X_{\sa}$.
That $\rho_{L,n}$ is finite follows from the norm compactness of 
$\mathcal{D}_1 (L)\cap\mathcal{B}^X_r$ for any $r>0$ 
(a consequence of \cite[Thm.\ 4.5]{GHDQMS} by scaling)  
along with Proposition~\ref{P-normLip-norm} below.

\begin{proposition}\label{P-completediambound}
The diameters of $UCP_n (X)$ relative to the respective metrics $\rho_{L,n}$ 
are finite and coincide for all $n\in\mathbb{N}$. 
\end{proposition}

\begin{proof}
The restriction map from $S(X)$ onto $S(\mathcal{D}(L))$ is evidently
a weak$^*$ homeomorphism which is isometric with respect to
$\rho_{L,1}$ and $\rho_L$ (Definition~\ref{D-cqms}), and so the
diameter of $S(X)$ with respect to $\rho_{L,1}$ is finite by
\cite[Thm.\ 4.5]{GHDQMS}. Now given $n\in\mathbb{N}$, $\varphi ,
\psi\in UCP_n (X)$, and $x\in X_{\sa}$ we can find a state $\sigma$
on $M_n$ such that
$$ | (\sigma\circ\varphi )(x) - (\sigma\circ\psi )(x) | = \| \varphi (x) -
\psi (x) \| . $$
It follows that the diameter of $UCP_n (X)$ is bounded above by that
of $S(X) = UCP_1 (X)$. On the other hand $S(X)$ embeds into $UCP_n (X)$
via the map which takes $\sigma\in S(X)$ to $x\mapsto\sigma (x)
1_{M_n}$, from which we see that the diameter of $UCP_n (X)$ is at least
that of $S(X)$, so that the two are equal.
\end{proof}

\begin{definition}\label{D-diameter}
Given a Lip-normed operator system $(X,L)$ we define its {\em diameter
$\diam (X,L)$} to be the common value of the diameters of $UCP_n (X)$
with respect to $\rho_{L,n}$ for $n\in\mathbb{N}$.
\end{definition} 

The next proposition, in addition to showing the finiteness of the
metrics $\rho_{L,n}$ (see the paragraph following 
Definition~\ref{D-ucpmetric}), will also be of use in 
Sections~\ref{S-distance} and \ref{S-tb} since it will enable us to 
streamline the statement and verification of conditions involving local 
approximation of elements of bounded Lip-norm. 

\begin{proposition}\label{P-normLip-norm}
Let $(X,L)$ be a Lip-normed operator system. Let $x\in\mathcal{D}(L)$ 
and let $r$ be the infimum of its spectrum.
Then $\| x - r 1 \| \leq L(x)\diam (X,L)$.
\end{proposition}

\begin{proof}
We can find $\sigma , \omega\in S(X)$ such that $\sigma (x - r 1) = 
\| x - r 1 \|$ and $\omega (x) = r$, whence 
\begin{align*}
\| x - r 1 \| & =  | \sigma (x-r 1) - \omega (x-r 1) | \\
&= | \sigma (x) - \omega (x) | \\
&\leq L(x) \diam (X,L) .
\end{align*}
\end{proof}

Since one of the requirements for a Lip-norm is that the associated
metric on the state space give rise to the weak$^*$ topology, one would
hope that the associated metrics on the matrix state spaces give rise to the
respective point-norm topologies. The following result shows that this is 
indeed the case. 

\begin{proposition}\label{P-ntopology}
The metric $\rho_{L,n}$ gives rise to the point-norm topology on $UCP_n (X)$.
\end{proposition}

\begin{proof}
Let 
$$ U_{\varphi , \Omega , \epsilon} = \{ \psi\in UCP_n (X) : \| \varphi (x) -
\psi (x) \| < \epsilon\text{ for all }x\in\Omega \} $$
be a basic open set in the point-norm topology, with $\varphi\in UCP_n (X)$,
$\epsilon > 0$, and $\Omega$ a finite subset of $A$. For each $x\in\Omega$
pick $y_{x,1} , y_{x,2}\in\mathcal{D}(L)$ with $\| y_{x,1} - \re (x) \| < 
\epsilon /6$ and $\| y_{x,2} - \im (x) \| < \epsilon /6$,
and choose $M>0$ so that $M\geq\max_{x\in\Omega ,j=1,2}L(y_{x,j})$. Now if
$\psi\in UCP_n (X)$ and 
$\rho_{L,n}(\varphi , \psi ) < (6M)^{-1}\epsilon$ then 
$\| \varphi (y_{x,j} ) - \psi (y_{x,j} ) \| < \epsilon /6$ for all $x\in\Omega$
and $j=1,2$, and hence
\begin{align*}
\| \varphi (\re (x)) - \psi ( \re (x)) \| &\leq 
\| \varphi (\re (x)) - \varphi (y_{x,1} ) \| + \| \varphi (y_{x,1} ) - 
\psi (y_{x,1} ) \| \\
&\hspace*{3cm}\ + \| \psi (y_{x,1} ) - \psi (\re (x)) \| \\
&< \epsilon /2 
\end{align*}
and similarly $\| \varphi (\im (x)) - \psi ( \im (x)) \| < \epsilon /2$,
so that $\| \varphi (x) - \psi (x) \| < \epsilon$.
Thus $U_{\varphi , \Omega , \epsilon}$ contains the open $\rho_{L,n}$-ball 
centred at $\varphi$ with radius $(6M)^{-1}\epsilon$, 
from which it follows that the metric topology is finer than the point-norm 
topology.

Suppose now that $\mathcal{B}(\varphi , \epsilon )$ is the $\rho_{L,n}$-ball
centred at some $\varphi\in UCP_n (X)$ with radius some $\epsilon > 0$.
Note that $\mathcal{D}_1 (L) \cap\mathcal{B}^X_{\diam (S(X))}$ is compact,
since $\mathcal{D}_1 (L)$ is closed by the definition of a Lip-normed
operator system and by \cite[Thm.\ 4.5]{GHDQMS} 
$\mathcal{D}_1 (L) \cap\mathcal{B}^X_1$
is totally bounded, which implies the total boundedness
of $\mathcal{D}_1 (L) \cap\mathcal{B}^X_{\diam (S(X))}$ via a scaling
argument. Hence we can find a 
finite set $\Omega\subset\mathcal{D}_1 (L) \cap\mathcal{B}^X_{\diam (S(X))}$ 
which is $(\epsilon /3)$-dense in $\mathcal{D}_1 (L) 
\cap\mathcal{B}^X_{\diam (S(X))}$. Thus if $\psi\in UCP_n (X)$
and $\| \varphi (x) - \psi (x) \| < \epsilon /3$ for all $x\in\Omega$ then
$\| \varphi (x) - \psi (x) \| < \epsilon$ for all $x\in\mathcal{D}_1 (L) \cap 
\mathcal{B}^X_{\diam (S(X))}$, and so $\mathcal{B}(\varphi , \epsilon )$ 
contains the point-norm basic open set
$$ \{ \psi\in UCP_n (X) : \| \varphi (x) -
\psi (x) \| < \epsilon /3 \text{ for all }x\in\Omega \} . $$
We conclude that the metric and point-norm topologies coincide on $UCP_n (X)$.
\end{proof}

We round out this section by showing that matrix state spaces embed 
isometrically under quotient maps, as do state spaces in the quantum metric 
setting. This will be crucial for the application 
of Hausdorff distance in formulating our matrix version of quantum 
Gromov-Hausdorff distance. 

\begin{proposition}\label{P-isometry}
Let $(X,L)$ be a Lip-normed operator system, $n\in\mathbb{N}$, 
$\Phi : X \to Y$ a surjective unital
$n$-positive map onto an operator system $Y$, and $L_Y$ the quotient Lip-norm 
on $Y$ induced by $L$ via $\Phi$. Then
the map $\Gamma : UCP_n (Y) \to UCP_n (X)$ given by
$\Gamma (\varphi ) = \varphi\circ\Phi$ is an isometry with respect
to $\rho_{L,n}$ and $\rho_{L_Y ,n}$.
\end{proposition}

\begin{proof}
Suppose $\Phi$ is $n$-positive and let $\varphi , \psi\in UCP_n (Y)$. 
Since $\Phi$ is Lip-norm-decreasing, we have $\rho_{L_Y ,n}(\varphi , \psi ) 
\geq\rho_{L,n}(\varphi\circ\Phi , \psi\circ\Phi )$. For the reverse
inequality, let $\epsilon > 0$ and choose $y\in\mathcal{D}_1 (L_Y )$
such that
$$ \rho_{L_Y ,n}(\varphi , \psi ) < \| \varphi (y) - \psi (y) \| + \epsilon 
. $$
We may assume $L(y) < 1$ for otherwise we can replace $y$ with
$\mu y$ for some $\mu < 1$ sufficiently close to $1$. Then by
definition of the quotient Lip-norm there is an $x\in\mathcal{D}_1
(L)$ such that $\Phi (x) = y$, and so
\begin{align*}
\rho_{L_Y ,n}(\varphi , \psi ) &< \| (\varphi\circ\Phi )(x) - (\psi\circ\Phi )
(x) \| + \epsilon \\
&\leq \rho_{L,n}(\varphi\circ\Phi , \psi\circ\Phi ) + \epsilon .
\end{align*}
Since $\epsilon$ was arbitrary, we conclude that
$\rho_{L_Y ,n}(\varphi , \psi ) = \rho_{L,n}(\varphi\circ\Phi , \psi\circ
\Phi )$, so that $\Gamma$ is an isometry with respect
to $\rho_{L_Y ,n}$ and $\rho_{L,n}$.
\end{proof}

\section{$n$-distance and complete distance}\label{S-distance}

The definition of quantum Gromov-Hausdorff distance \cite{GHDQMS} involves
forming a direct sum and considering Lip-norms thereupon which induce
the given Lip-norms on the summands. One then takes an infimum of the 
Hausdorff distances between the state spaces under their isometric 
embeddings into the state space of the direct sum. We will apply the same 
procedure here with regard to the matrix state spaces.

\begin{notation}\label{N-induce}
Let $(X,L_X )$ and $(Y,L_Y )$ be Lip-normed operator
systems. We denote by $\mathcal{M}(L_X , L_Y )$
the collection of closed Lip-norms on 
$\mathcal{D}(L_X ) \oplus\mathcal{D}(L_Y )$ 
which induce $L_X$ and $L_Y$ via the quotient maps onto
$\mathcal{D}(L_X )$ and $\mathcal{D}(L_Y )$, respectively.
\end{notation}

Let $(X,L_X )$ and $(Y,L_Y )$ be Lip-normed operator systems
and $L\in\mathcal{M}(L_X , L_Y )$. Since the projection map 
$X\oplus Y\to X$ is u.c.p., by Proposition~\ref{P-isometry} we obtain an 
isometry $UCP_n (X)\to UCP_n (X\oplus Y)$ with respect to $\rho_{L_X}$
and $\rho_L$. Similarly, we also have an isometry $UCP_n (Y)\to 
UCP_n (X\oplus Y)$. For notational simplicity we will thus identify 
$UCP_n (X)$ and $UCP_n (Y)$ with their respective images under these 
isometries.

\begin{definition}\label{D-dist}
Let $(X,L_X )$ and $(Y,L_Y )$ be Lip-normed operator systems.
For each $n\in\mathbb{N}$ we define the {\em $n$-distance}
$$ \dist^n_s (X,Y) = \inf_{L\in\mathcal{M}(L_X , L_Y 
)} \dist^{\rho_{L,n}}_H (UCP_n (X) , UCP_n (Y)) $$
where $\dist^{\rho_{L,n}}_H$ denotes Hausdorff distance with respect to
the metric $\rho_{L,n}$. We also define the {\em complete distance}
$$ \dist_s (X,Y) = \inf_{L\in\mathcal{M}(L_X , L_Y )} \,\,
\sup_{n\in\mathbb{N}}\, 
\dist^{\rho_{L,n}}_H (UCP_n (X) , UCP_n (Y)) . $$
\end{definition}

\begin{proposition}[triangle inequality]\label{P-triangle}
If $(X,L_X )$, $(Y,L_Y )$, and $(Z,L_Z )$ are Lip-normed operator 
systems then
$$ \dist^n_s (X,Z) \leq \dist^n_s (X,Y) + \dist^n_s (Y,Z) , $$
for all $n\in\mathbb{N}$, and
$$ \dist_s (X,Z) \leq \dist_s (X,Y) + \dist_s (Y,Z) . $$
\end{proposition}

\begin{proof}
This follows by exactly
the same argument used for quantum Gromov-Hausdorff distance in 
\cite{GHDQMS}, since in the last part of the proof of 
\cite[Thm.\ 4.3]{GHDQMS} we can replace the state spaces by matrix state
spaces and the reference to \cite[Prop.\ 3.1]{GHDQMS} by a reference to
our Proposition~\ref{P-isometry}. 
\end{proof} 

In order to build a general framework for estimating distance between 
quantum metric spaces, Rieffel formulates in \cite[Defn.\ 5.1]{GHDQMS} the
notion of a {\em bridge}, which we now recall.
 
\begin{definition}\label{D-bridge}
Let $(A, L_A )$ and $(B, L_B )$ be compact quantum metric spaces. A
{\em bridge} between $(A, L_A )$ and $(B, L_B )$ is a norm-continuous
semi-norm $N$ on $A\oplus B$ such that 
\begin{enumerate}
\item[(i)] $N(1_A , 1_B ) = 0$ while $N(1_A , 0) \neq 0$, and
\item[(ii)] for each $a\in A$ and $\delta > 0$ there exists a $b\in B$
such that 
$$ \max (L_B (b) , N(a,b )) \leq L_A (a) + \delta , $$ 
with the same statement also holding upon interchanging $A$ and $B$.
\end{enumerate}
\end{definition}

Theorem 5.2 in \cite{GHDQMS} then shows that if $N$ is a bridge between the 
compact quantum metric spaces $(A, L_A )$ and $(B, L_B )$ then 
$$ L(a,b) = \max (L_A (a) , L_B (b) , N(a,b)) $$
defines a Lip-norm $L$ on $A\oplus B$ which induces $L_A$ and $L_B$ via 
the respective quotient maps. Since $N$ is norm-continuous, $L$ will
be closed if $L_A$ and $L_B$ are both closed.

If $(X,L_X )$ and $(Y, L_Y )$ are Lip-normed
operator systems then by a {\em bridge between 
$(X,L_X )$ and $(Y,L_Y )$} we will mean a bridge between
the compact quantum metric spaces $(\mathcal{D}(L_X ), L_X )$ and
$(\mathcal{D}(L_Y ), L_Y )$. We begin by illustrating this notion with a
simple example which shows that if we scale a Lip-norm by a factor $\lambda$
and let $\lambda\to\infty$ then we obtain convergence to a ``point,'' just
as for ordinary metric spaces.

\begin{example}
Let $(X,L)$ be a Lip-normed operator system. For each $\lambda > 0$
define the Lip-normed operator system $(X,L_\lambda )$ by setting
$L_\lambda = \lambda L$. Let $(\mathbb{C} , P)$ be the ``one-point''
Lip-normed operator system, with $P(\mu )=0$ for all $\mu\in
\mathbb{C}$. Then
$$ \dist_s ((X,L),(\mathbb{C} , P)) \leq C\lambda^{-1} $$
where $C=\diam (X,L)$.
To show this we define a bridge on $\mathcal{D}(L_X ) \oplus\mathbb{C}$ by
$$ N_\lambda (x,\mu ) = C^{-1}\lambda \| x - \mu 1_X \| . $$
To see that this is indeed a bridge we verify condition (ii) in 
Definition~\ref{D-bridge} by observing that $N_\lambda (\mu 1_X , \mu ) = 0$ 
for all $\mu\in\mathbb{C}$ while
if $x\in\mathcal{D}(L_X )$ then letting $r$ denote the infimum of
the spectrum of $x$ we have by Proposition~\ref{P-normLip-norm}
$$ \| x - r 1_X \| \leq \diam (X,L) L_X (x) = C\lambda^{-1}L_\lambda (x) . $$
Let $M_\lambda$ be the Lip-norm in $\mathcal{M}(L_\lambda , P)$ given by
$M_\lambda (x,\mu ) = \max (L_\lambda (x) , P(\mu ) , N_\lambda (x,\mu ))$.
Let $n\in\mathbb{N}$ and let $\psi$ be the unique element in 
$UCP_n (\mathbb{C})$. If $\varphi\in UCP_n (X)$ and $(x,\mu )\in
\mathcal{D}_1 (M)$ then
$$ \| \varphi (x) - \psi (\mu ) \| = \| \varphi (x-\mu 1_X ) \| \leq
\| x-\mu 1_X \| \leq C\lambda^{-1} , $$
yielding the desired complete distance estimate. Hence 
$(X,L_\lambda )$ converges to $(\mathbb{C} , P)$ as $\lambda\to\infty$ for
complete distance.
\end{example}

As in the quantum metric setting \cite[Prop.\ 5.4]{GHDQMS} we can apply 
the concept of a bridge to show that the complete distance (and hence also 
the $n$-distance) is always finite.

\begin{proposition}\label{P-distdiambound}
If $(X,L_X )$ and $(Y,L_Y )$ are Lip-normed operator systems then
$$ \dist_s (X,Y) \leq\diam (X,L_X ) + \diam (Y,L_Y ) . $$
\end{proposition}

\begin{proof}
As in the proof of \cite[Prop.\ 5.4]{GHDQMS}, for arbitrary 
$\gamma > 0$, $\sigma_0 \in S(X)$, and $\omega_0 \in S(Y)$ we can 
construct a bridge 
$$ N(x,y) = \gamma^{-1} | \sigma_0 (x) - \omega_0 (y) | . $$
Let $L$ be the Lip-norm in $\mathcal{M}(L_X , L_Y )$ given by
$L(x,y) = \max (L_X (x) , L_Y (y) , N(x,y))$.
Then if $(x,y)\in\mathcal{D}_1 (L)$, $\varphi\in UCP_n (X)$, and
$\psi\in UCP_n (Y)$ we can find a $\sigma\in S(M_n )$ such that
$\| \varphi (x) - \psi (y) \| = | (\sigma\circ\varphi )(x) -
(\sigma\circ\psi )(x) |$ whence
\begin{align*}
\| \varphi (x) - \psi (y) \| &\leq | (\sigma\circ\varphi )(x) -
\sigma_0 (x) | + | \sigma_0 (x) - \omega_0 (y) | + | \omega_0 (y) -
(\sigma\circ\psi )(y) | \\
&\leq \rho_{L_X ,1} (\sigma\circ\varphi , \sigma_0 ) + \gamma +
\rho_{L_Y ,1} (\omega_0 , \sigma\circ\psi ) .
\end{align*}
Since $\gamma$ was arbitrary the proposition follows.
\end{proof}

Propositions~\ref{P-distest} and \ref{P-quotientest}
yield estimates on the complete distance in situations involving bridges 
constructed via the norm. 

\begin{proposition}\label{P-distest}
Let $(X,L_X )$ and $(Y,L_Y )$ be Lip-normed operator systems,
and suppose $X$ and $Y$ are operator subsystems of an operator
system $Z$. Let $\epsilon >0$ and
suppose that $N(x,y) = \epsilon^{-1}\| x-y \|$ defines a bridge 
between $(X,L_X )$ and $(Y,L_Y )$. Then $\dist_s (X,Y) \leq\epsilon$.
\end{proposition}

\begin{proof}
Let $L$ be the Lip-norm in $\mathcal{M}(L_X , L_Y )$ given by
$$ L(x,y) = \max\left( L_X (x) , L_Y (y) , \epsilon^{-1}
\| x - y \| \right) $$ 
for all $(x,y)\in \mathcal{D}(L_X )\oplus\mathcal{D}(L_Y )$.
Let $n\in\mathbb{N}$ and $\varphi\in UCP_n (X)$. By Arveson's extension
theorem we can extend $\varphi$ to a u.c.p.\ map $\tilde{\varphi} : Z\to
M_n$. Then if $(x,y) \in\mathcal{D}_1 (L)$ we have 
$$ \| \varphi (x) - \tilde{\varphi}(y) \| \leq\| x - y \| < \epsilon $$
and thus $\rho_{L,n}(\varphi , \tilde{\varphi}|_Y ) < \epsilon$. We can
interchange the roles of $X$ and $Y$ and apply the same argument 
to conclude that $\dist_s (X,Y) \leq\epsilon$.
\end{proof}

\begin{proposition}\label{P-quotientest}
Let $(X,L_X )$ be a Lip-normed operator system, $Y$ an operator system,
$\Phi : X\to Y$ a surjective unital positive map, and $L_Y$ the
quotient Lip-norm induced by $L_X$ via $\Phi$. Let $Z$ be an operator 
system containing $X$ as an operator subsystem and let $\Gamma : Y\to Z$ 
be a unital map such that $\| (\Gamma\circ\Phi )(x) - x \| \leq\epsilon$ 
for all $x\in\mathcal{D}_1 (L)$. If $\Phi$ and $\Gamma$ are $n$-positive 
then $\dist^n_s (X,Y)\leq\epsilon$, and if $\Phi$ and $\Gamma$ are 
completely positive then $\dist_s (X,Y)\leq\epsilon$.
\end{proposition}

\begin{proof}
Given $\eta > 0$ we define the bridge $N$ between $(X,L)$ and $(Y,L_Y )$ by
$N(x,y) = \eta^{-1}\| \Phi (x) - y \|$
for all $(x,y)\in X\oplus Y$ (to verify condition (ii) of 
Definition~\ref{D-bridge}, given $x\in X$ and $\delta > 0$ we can choose 
$y=\Phi (x)$, while if $y\in Y$ and $\delta > 0$ we can take any $x\in X$ 
such that $\Phi (x) = y$ and $L_X (x) \leq L_Y (y) + \delta$).
Let $L$ be the Lip-norm $L$ in $\mathcal{M}(L_X , L_Y )$ given by 
$$ L(x,y) = \max (L_X (x) , L_Y (y) , N(x,y)) $$
for all $(x,y)\in X\oplus Y$. Suppose that $\Phi$ is $n$-positive.
Then if $\varphi\in UCP_n (Y)$ we have $\varphi\circ\Phi\in UCP_n (X)$, 
and so if $L(x,y) \leq 1$ then $\| \Phi (x) - y \| \leq\eta$ so that 
$\| (\varphi\circ\Phi )(x) - \varphi (y) \| \leq\eta$, whence
$\rho_{L,n}(\varphi\circ\Phi , \varphi ) \leq\eta$.
On the other hand if $\varphi\in UCP_n (X)$ then, extending $\varphi$ to 
a u.c.p.\ map $\tilde{\varphi} : Z\to M_n$ by Arveson's extension theorem, 
we have $\tilde{\varphi}\circ\Gamma\in UCP_n (Y)$, and so if $L(x,y) \leq 1$ 
then 
\begin{align*}
\| \varphi (x) - (\tilde{\varphi}\circ\Gamma )(y) \| &\leq \| \tilde{\varphi} 
(x - (\Gamma\circ\Phi )(x)) \| + \| (\tilde{\varphi}\circ\Gamma )
(\Phi (x) - y) \| \\
&\leq \epsilon + \eta ,
\end{align*}
yielding $\rho_{L,n}(\varphi , \tilde{\varphi}\circ\Gamma ) \leq\epsilon + 
\eta$. Since $\eta$ was arbitrary we conclude that 
$\dist^n_s (X,Y)\leq\epsilon$. In the case that $\Phi$ is completely 
positive we can apply the above argument over all $n\in\mathbb{N}$ to obtain 
$\dist_s (X,Y)\leq\epsilon$.
\end{proof}

The following three propositions guarantee approximability by Lip-normed
well-behaved finite-dimensional operator systems under conditions on 
the given operator system or $C^*$-algebra which hold in a wide range of
situations.
 
\begin{proposition}\label{P-finiterank}
Let $(X,L_X )$ be a Lip-normed nuclear operator system. Then 
for every $\epsilon > 0$ there is a Lip-normed 
operator system $(Y,L_Y )$ such that $Y$ is an operator subsystem of a 
matrix algebra and $\dist_s (X,Y) \leq\epsilon$.
\end{proposition} 

\begin{proof}
Let $\epsilon > 0$. Since $X$ is nuclear and the set $\mathcal{D}_1 
(L_X ) \cap X_{\diam (X,L_X )}$ is compact (see the second half of the proof 
of Proposition~\ref{P-ntopology}), 
by Proposition~\ref{P-normLip-norm} we can find a matrix algebra $M_k$ and
u.c.p.\ maps $\Phi : X \to M_k$ and $\Gamma : M_k \to X$ 
such that $\| (\Gamma\circ\Phi )(x) - x \| \leq\epsilon$
for all $x\in\mathcal{D}_1 (L_X )$. Consider the image $Y$ of $\Phi$ 
and the resulting quotient Lip-norm $L_Y$ on $Y$. Then
by Proposition~\ref{P-quotientest} we have $\dist_s (X,Y) \leq\epsilon$, 
yielding the result.
\end{proof}

By a proof similar to that of Proposition~\ref{P-finiterank}, we also
have the following.

\begin{proposition}\label{P-exact}
Let $(A,L)$ be a Lip-normed unital exact $C^*$-algebra. Then 
for every $\epsilon > 0$ there is a Lip-normed 
operator system $(Y,L_Y )$ such that $Y$ is an operator subsystem of a 
matrix algebra and $\dist_s (A,Y) \leq\epsilon$.
\end{proposition}

A separable $C^*$-algebra $A$ is said to be a {\em strong NF algebra} if
it is the inductive limit of a generalized inductive system $(A_n ,
\phi_{n,m} )$ with each $A_n$ a finite-dimensional $C^*$-algebra and
each $\phi_{n,m}$ a complete order embedding \cite[Defn.\ 5.2.1]{BK}
(a {\em complete order embedding} from a $C^*$-algebra $B$ to a
$C^*$-algebra $A$ is a c.p.\ isometry $\Phi : A\to B$ such that 
$\Phi^{-1} : \Phi (A) \to B$ is a c.p.\ map). 

\begin{proposition}\label{P-strongNF}
Let $(A,L)$ be a Lip-normed unital strong NF algebra. Then for every 
$\epsilon > 0$ there is a Lip-normed finite-dimensional $C^*$-algebra 
$(B,L_B )$ such that $\dist_s (A,B) \leq\epsilon$.
\end{proposition}

\begin{proof}
Since the set $\mathcal{D}_1 (L) \cap X_{\diam (A,L)}$ is
compact (see the second half of the proof of Proposition~\ref{P-ntopology})
and $A$ is strong NF, by \cite[Thm.\ 6.1.1]{BK} and
Proposition~\ref{P-normLip-norm} we can find a finite-dimensional
$C^*$-algebra $B$, a unital complete order embedding
$\Gamma : B \to A$, and a (surjective) u.c.p.\ map $\Phi : A \to B$ such 
that $\Phi\circ\Gamma = \id_B$ and $\| (\Gamma\circ\Phi )(a) - a \| \leq
\epsilon$ for all $a\in\mathcal{D}_1 (L)$. Then $L$ induces a Lip-norm 
$L_B$ on $B$ via $\Phi$, and $\dist_s (A,B) \leq
\epsilon$ by Proposition~\ref{P-quotientest}.
\end{proof}

Proposition~\ref{P-strongNF} applies for instance to noncommutative
tori Lip-normed via the ergodic action of ordinary tori, as described in
Example~\ref{E-tori}. In this situation, however, one would hope to be able 
to approximate by finite-dimensional $C^*$-algebras Lip-normed 
via models of the original action, as in the following example. 

\begin{example}\label{E-sphere}
In \cite{MACS} Rieffel shows, in the context of Berezin quantization,
that the sphere $S^2$ is a limit of matrix algebras with 
respect to quantum Gromov-Hausdorff distance.
In fact a more general statement applying to integral coadjoint orbits
of a compact Lie group is established. We will briefly indicate how
Rieffel's approach leads to precisely the same estimates for complete 
distance, adopting the same notation as in \cite{MACS}, to which we
refer the reader for more details. 

Given a compact group $G$ consider the $C^*$-algebra $B$ of all bounded
operators on a Hilbert space on which $G$ is irreducibly and unitarily
represented. Given a rank-one projection $P\in B$ we define for each
$T\in B$ the Berezin covariant symbol $\sigma_T$ with respect to $P$ by
$\sigma_T (x) = \tau (T\alpha_x (P))$ where $\tau$ is the unnormalized
trace on $B$ and $\alpha$ is the action of $G$ on $B$ given by conjugation.
Denoting by $H$ the stability subgroup of $P$ for $\alpha$, we thereby 
obtain a map $\sigma$ from $B$ to $A=C(G/H)$ which is unital and positive,
and hence u.c.p.\ since the range is a commutative $C^*$-algebra. The
action $\alpha$ along with a length function $\ell$ give rise to a
Lip-norm $L_B$ on $B$ as in Example~\ref{E-actions}, and similarly
the action of $G$ on $G/H$ by left translation combines with $\ell$
to produce a Lip-norm $L_A$ on $A$ (permitting the value $+\infty$ for
convenience). Corollary 2.4 of \cite{MACS} shows that, for some
$\gamma > 0$, there is a bridge between $(A,L_A )$ and $(B,L_B )$ of the
form
$$ N(f,T) = \gamma^{-1} \| f - \sigma_T \|_{\infty} . $$
Proposition 1.3 of \cite{MACS} then shows that $S(A)$ is in the
$\gamma$-neighbourhood of $S(B)$ under the metric defined by the Lip-norm
$L$ on $A\oplus B$ associated to $N$. But the argument there also applies
at the matrix level: given $\varphi\in UCP_n (A)$ we have
$\varphi\circ\sigma\in UCP_n (B)$ since $\sigma$ is u.c.p., 
and if $L(f,T)\leq 1$ then $\| f - \sigma_T \|_{\infty} \leq\gamma$ so that
$$ \| \varphi (f) - (\varphi\circ\sigma )(T) \| = \| \varphi (f-\sigma_T ) \|
\leq \| f - \sigma_T \|_{\infty} \leq\gamma , $$
showing that $UCP_n (A)$ lies in the $\gamma$-neighbourhood of $UCP_n (B)$
with respect to the metric $\rho_{L,n}$. Now on the other hand if $\psi\in
UCP_n (B)$ then we can consider the adjoint operator $\breve{\sigma} :
A\to B$ and take the composition $\psi\circ\breve{\sigma}$, 
which is in $UCP_n (A)$ since $\breve{\sigma}$ is unital and positive 
and hence u.c.p.\ because its
domain is a commutative $C^*$-algebra. Then if $L(f,T)\leq 1$ we have 
$\| f - \sigma_T \|_{\infty} \leq\gamma$ from which we obtain the estimate
\begin{align*}
\| (\psi\circ\breve{\sigma})(f) - \psi (T) \| &= \| \psi (\breve{\sigma}_f 
- T ) \| \\
&\leq \| \breve{\sigma}_f - T \| \\
&\leq \| \breve{\sigma}_f - \breve{\sigma}(\sigma_T ) \| + 
\| \breve{\sigma}(\sigma_T ) - T \| \\
&\leq \| f - \sigma_T \|_{\infty} + \| \breve{\sigma}(\sigma_T ) - T \| \\
&\leq \gamma + \| \breve{\sigma}(\sigma_T ) - T \| 
\end{align*}
exactly as in the case $n=1$ in the discussion following
\cite[Cor.\ 2.4]{MACS}. When $G$ is a compact Lie group Rieffel shows
that, by replacing $B$ with the $C^*$-algebra of bounded 
operators on the $m$th tensor power of the original Hilbert space,
both the bridge constant $\gamma$ and the term $\| \breve{\sigma}(\sigma_T ) 
- T \|$ can be made arbitrarily small by taking $m$ sufficiently large
(see Theorem 3.2 and Sections 3--5 of \cite{MACS}), yielding
an asymptotically vanishing bound on the quantum Gromov-Hausdorff distance
as a result of the estimates in the two previous displays 
for the case $n=1$. But since these estimates apply equally well for all
$n$ we get the same bounds for complete distance.
\end{example}

\section{Distance zero}\label{S-zero}

This section is aimed at establishing that $\dist^n_s (X,Y) = 0$
(resp.\ $\dist_s (X,Y) = 0$) is equivalent to the existence of a 
bi-Lip-isometric unital $n$-order isomorphism (resp.\ bi-Lip-isometric unital 
complete order isomorphism) between $X$ and $Y$ (Theorem~\ref{T-orderisom}). 
One direction is straightforward:

\begin{proposition}\label{P-Lipisometric}
Let $(X,L_X )$ and $(Y,L_Y )$ be Lip-normed operator systems. If there
is a bi-Lip-isometric unital $n$-order isomorphism between $X$ and $Y$ then
$\dist^n_s (X,Y) = 0$. If there is a bi-Lip-isometric unital complete
order isomorphism between $X$ and $Y$ then $\dist_s (X,Y) = 0$.
\end{proposition}

\begin{proof}
This follows from Proposition~\ref{P-quotientest}, taking $\Phi$ there to 
be a bi-Lip-isometric unital $n$-order 
isomorphism (resp.\ bi-Lip-isometric unital complete order isomorphism) 
from $X$ onto $Y$ and taking $\Gamma$ to be its inverse. 
\end{proof}

For the converse, it will be convenient to extend our Lip-norms in an
adjoint-invariant way (Definition~\ref{D-semi-norm}) and to introduce
a collection of matrix semi-norms (Definition~\ref{D-matrixsn}).

\begin{definition}\label{D-semi-norm}
Let $(X,L)$ be a Lip-normed operator system. 
We define the semi-norm $L^e$ on $X$ by
$$ L^e (x) = \sup \left\{ \frac{| \sigma (x) - \omega (x) 
|}{\rho_{L,1}(\sigma , \omega )} : \sigma , \omega\in S(X) \text{ and }
\sigma\neq\omega \right\} $$
for all $x\in X$ (permitting the value $+\infty$). 
\end{definition}

\begin{proposition}\label{P-restriction}
The set of self-adjoint elements on which the semi-norm $L^e$ in
Definition~\ref{D-semi-norm} is finite coincides
with the domain of $L$, and on this set $L^e = L$.
\end{proposition}

\begin{proof}
By Proposition 6.1 of \cite{GHDQMS}, $\mathcal{D}(L)$ corresponds to the 
subspace of affine functions on $S(\mathcal{D}(L))$ which are Lipschitz for
$\rho_L$. This, along with the fact that
the restriction map from $S(X)$ to $S(\mathcal{D}(L))$ is a weak$^*$
homeomorphism which is isometric for $\rho_{L,1}$ and $\rho_L$ (see
the proof of Proposition~\ref{P-selfadj}), implies the result.
\end{proof}

\begin{definition}\label{D-matrixsn}
Let $(X,L)$ be a Lip-normed operator system and $n\in\mathbb{N}$. We define
the semi-norm $L^n$ on $M_n \otimes X$ by 
$$ L^n (x) = \max_{1\leq i,j \leq n} L^e (x_{ij}) $$
for all $x = (x_{ij}) \in M_n (X) \cong M_n \otimes X$.
\end{definition}

For the meaning of the notation $\mathcal{D} (\cdot )$ and 
$\mathcal{D}_\lambda (\cdot )$, as will be applied to the semi-norms $L^e$ 
and $L^n$, see Notation~\ref{N-L}.

\begin{lemma}\label{L-dense}
Let $(X,L)$ be a Lip-normed operator system. Then $\mathcal{D}(L^e )$ is
dense in $X$ and $\mathcal{D}(L^n )$ is dense in $M_n \otimes X$ for
all $n\in\mathbb{N}$.
\end{lemma}

\begin{proof}
Since $L^e$ is adjoint-invariant and coincides with $L$ on the dense 
real subspace $\mathcal{D}(L)$ of $X_{\sa}$ by 
Proposition~\ref{P-restriction}, using the decomposition of elements into 
real and imaginary parts we see that $\mathcal{D}(L^e )$ is dense in $X$.
As a direct consequence $\mathcal{D}(L^n )$ is dense in $M_n \otimes X$ for
all $n\in\mathbb{N}$.
\end{proof}

\begin{lemma}\label{L-nbound}
Let $(X,L_X )$ and $(Y,L_Y )$ be Lip-normed operator systems, and
let $L\in\mathcal{M}(L_X , L_Y )$ and $n\in\mathbb{N}$. Set
$r = \dist^{\rho_{L,n}}_H (UCP_n (X) , UCP_n (Y))$. Then for 
every $\psi\in SCP_n (Y)$ there is a $\varphi\in SCP_n (X)$ such that,
for all $(x,y)\in\mathcal{D}(L^e )$,
$$ \| \varphi (x) - \psi (y) \| \leq 2n^3 L^e (x,y) r . $$
\end{lemma}

\begin{proof} 
To prove the lemma we may assume that $\psi (1)$ has 
full support in $M_n$, for otherwise for every $k\in\mathbb{N}$ we can 
perturb $\psi$ to a convex combination 
$(1-k^{-1})\psi + k^{-1}\alpha$ where $\alpha (x) = \omega (x) 
1_{M_n}$ for some $\omega\in S(Y)$ and all $x\in X$ (in which case the 
corresponding state $\sigma_\psi$ on $M_n \otimes Y$ is perturbed to 
another a state), find a suitable $\varphi_k$ as in the lemma statement
with respect to $\psi_k$, and then take a point-norm limit point of 
$\{ \varphi_k \}_{k\in\mathbb{N}}$ to obtain the desired $\phi$. We can 
thus consider the map $\psi' \in UCP_n (Y)$ given by
$$ \psi' (y) = \psi (1)^{-\frac12} \psi (y) \psi (1)^{-\frac12} . $$
for all $y\in Y$. By assumption we can find a $\varphi' \in UCP_n (X)$
such that
$$ \rho_{L,n} (\varphi' , \psi' ) \leq r . $$
Let $\varphi : Y\to M_n$ be the c.p.\ map given by
$$ \varphi (x) = \psi (1)^{\frac12} \varphi' (x) \psi (1)^{\frac12} $$
for all $x\in X$. Then $\varphi (1) = \psi (1)$, which implies that 
$\sigma_\varphi$ is a state on $M_n \otimes X$, so that
$\varphi \in SCP_n (Y)$. Since $\sigma_\psi$ is a state on $M_n \otimes Y$ 
we must have $\| \psi (1) \| \leq n^3$, and 
thus if $(x,y)\in\mathcal{D}(L^e )$ then
\begin{align*}
\| \varphi (x) - \psi (y) \| &= \| \psi (1)^{\frac12} (\varphi' (x) - 
\psi' (y) ) \psi (1)^{\frac12} \| \\
&\leq \| \psi (1)^{\frac12} \| \| \varphi' (x) - \psi' (y) \| 
\| \psi (1)^{\frac12} \| \\
&\leq n^3 (\| \varphi' (\re (x)) - \psi' (\re (y)) \| + \|
\varphi' (\im (x)) - \psi' (\im (y)) \| ) \\
&\leq n^3 (L(\re (x), \re (y)) + L(\im (x) , \im (y)))
\rho_{L,n}(\varphi' , \psi' ) \\
&\leq 2n^3 L^e (x,y) \rho_{L,n}(\varphi' , \psi' ) \\
&\leq 2n^3 L^e (x,y) r
\end{align*}
with the second last inequality following from the adjoint invariance of
$L^e$ and the fact that $L^e = L$ on $\mathcal{D}(L)$ by 
Proposition~\ref{P-restriction}.
\end{proof}

\begin{definition}\label{D-N}
Let $(X,L_X )$ and $(Y,L_Y )$ be Lip-normed operator systems, and let 
$L\in\mathcal{M}(L_X , L_Y )$ and $n\in\mathbb{N}$. 
For each $\lambda\geq 0$, $x\in\mathcal{D}(L^n_X )$, and
$y\in\mathcal{D}(L^n_Y )$ we set
\begin{gather*}
\mathcal{N}^\lambda_{L^n ,Y}(x) = \left\{ z\in M_n \otimes Y : (x,z) \in
\mathcal{D}_\lambda (L^n ) \right\} , \\
\mathcal{N}^\lambda_{L^n ,X}(y) = \left\{ z\in M_n \otimes X : (z,y) \in
\mathcal{D}_\lambda (L^n ) \right\} , \\
\end{gather*}
\end{definition}

\begin{lemma}\label{L-diambound}
Let $(X,L_X )$ and $(Y,L_Y )$ be Lip-normed operator systems, and
let $L\in\mathcal{M}(L_X , L_Y )$ and $n\in\mathbb{N}$. Set
$r = \dist^{\rho_{L,n}}_H (UCP_n (X) , UCP_n (Y))$.
If $x\in\mathcal{D}(L^n_X )$ and $\lambda > 2L^n_X (x)$ then
\begin{enumerate}
\item[(i)] $\mathcal{N}^\lambda_{L^n ,Y}(x)$ is non-empty and closed,
\item[(ii)] the norms of elements of $\mathcal{N}^\lambda_{L^n ,Y}(x)$ 
are bounded by $4(\| x \| + \lambda n^4 r)$, and if $x$ and $y\in
\mathcal{N}^\lambda_{L^n ,Y}(x)$ are self-adjoint then the norm of $y$ is
bounded by $\| x \| + 2\lambda n^4 r$, 
\item[(iii)] the norm diameter of $\mathcal{N}^\lambda_{L^n ,Y}(x)$
is bounded by $8\lambda n^4 r$, 
\item[(iv)] if $x'$ is another element of $\mathcal{D}(L^n_X )$ and $\lambda$ 
is also strictly larger than $L^n_X (x' )$ then the Hausdorff
distance between $\mathcal{N}^\lambda_{L^n ,Y}(x)$ and
$\mathcal{N}^\lambda_{L^n ,Y}(x' )$ is bounded by $8\lambda n^4 r + 
4\| x-x' \|$, and 
\item[(v)] if $x\geq 0$ then there is a (self-adjoint) 
$y\in\mathcal{N}^\lambda_{L^n ,Y}(x)$ with $y \geq -2\lambda n^4 r$.
\end{enumerate}
Statements (i)--(v) also hold for $y\in\mathcal{D}(L^n_Y)$ and 
$\lambda > 2 L^n_Y (y)$ if $L_X$ and
$\mathcal{N}^\lambda_{L^n ,Y}(x)$ are
replaced by $L_Y$ and $\mathcal{N}^\lambda_{L^n ,X}(y)$, respectively.
\end{lemma}

\begin{proof}
Since symmetry will take care of the last sentence of the proposition
statement, we prove (i)--(v) as written for $x\in\mathcal{D}(L^n )$ and 
$\lambda > 2 L^n_X (x)$. We begin with (i). For $1\leq i,j \leq n$ we
have $\max (L_X (\re (x_{ij})) , L_X (\im (x_{ij})) < L^e (x_{ij}) <
\lambda /2$ so that we can find $y_{ij} \in\mathcal{N}^{\lambda /2}_{L^n ,Y}
(\re (x_{ij}))$ and $z_{ij} \in\mathcal{N}^{\lambda /2}_{L^n ,Y}
(\im (x_{ij}))$ since $L_X$ is the quotient Lip-norm induced by $L$. Then
\begin{align*}
L^n ((x_{ij}) , (y_{ij}) + i (z_{ij})) &\leq 
L^n ((\re (x_{ij})) , (y_{ij})) + L^n ((\im (x_{ij})) , 
(z_{ij})) \\
&< \lambda
\end{align*}
so that $\mathcal{N}^{\lambda /2}_{L^n ,Y}(x)$ contains 
$(y_{ij}) + i (z_{ij})$\
and is in particular non-empty. That this set is also closed follows from
the lower semicontinuity of $L^e$, which is easily checked.

For (ii), let $y=(y_{ij})\in\mathcal{N}^\lambda_{L^n ,Y}(x)$ and  
$\psi\in SCP_n (Y)$. By Lemma~\ref{L-nbound} there is a $\varphi\in SCP_n (X)$
such that $\| \varphi (z) - \psi (w) \| \leq 2n^3 L^e (z,w)r$ for all
$z,w\in\mathcal{D}(L^e )$. We then have
\begin{align*}
| \sigma_\varphi (x) - \sigma_\psi (y) | &= \frac1n \Bigg| \sum_{i,j}
(\varphi (x_{ij})_{ij} - \psi (y_{ij})_{ij}) \Bigg| \\
&\leq \frac1n \sum_{i,j} | \varphi (x_{ij})_{ij} - \psi (y_{ij})_{ij}| \\
&\leq \frac1n \sum_{i,j} \| \varphi (x_{ij}) - \psi (y_{ij}) \|  \\
&\leq 2 n^2 \sum_{i,j} L^e (x_{ij} , y_{ij}) r \\
&\leq 2 n^4 L^n (x,y) r \\
&\leq 2 \lambda n^4 r .
\end{align*}
It follows that $| \sigma_\psi (y) | \leq 2(\| x \| + \lambda n^4 r)$, 
and so $| \sigma_\psi (\re (y)) |$ and $| \sigma_\psi (\im (y)) |$ are both 
bounded by $2(\| x \| + \lambda n^4 r)$, from which we conclude that 
$$ \| y \| \leq \| \re (y) \| + \| \im (y) \| \leq 4(\| x \| + 
\lambda n^4 r) . $$
If $x$ and $y\in\mathcal{N}^\lambda_{L^n ,Y}(x)$ are self-adjoint 
then the above argument shows that the norm of $y$ is in fact bounded by 
$\| x \| + 2\lambda n^4 r$.

To establish (iii), suppose
$y , y' \in M_n \otimes Y$ are such that $(x,y) , (x,y' )\in
\mathcal{D}_\lambda (L^e )$. Let $\psi\in SCP_n (Y)$. Arguing as in the proof 
of (ii), there exists by Lemma~\ref{L-nbound} a $\varphi \in SCP_n (X)$ 
such that $\| \varphi (z) - \psi (w) \| \leq n^3 L^e (z,w)r$ for all
$z,w\in\mathcal{D}(L^e )$, so that both
$| \sigma_\varphi (x) - \sigma_\psi (y) |$ and
$| \sigma_\varphi (x) - \sigma_\psi (y' ) |$ are bounded by $2\lambda n^4 r$, 
whence $| \sigma_\psi (y-y' ) | \leq 4\lambda n^4 r$.
It follows that $\| y-y' \| \leq 8\lambda n^4 r$,
and so we obtain (iii).

For (iv), suppose $y\in\mathcal{N}^\lambda_{L^n ,Y}(x)$ and
$y' \in\mathcal{N}^\lambda_{L^n ,Y}(x' )$. Then as in the proof of (ii)
we can find a $\varphi\in SCP_n (X)$ such that both $| \sigma_\varphi (x) - 
\sigma_\psi (y) |$ and $| \sigma_\varphi (x' ) - \sigma_\psi (y' ) |$ are
bounded by $2\lambda n^4 r$, and the triangle inequality yields
\begin{align*}
| \sigma_\psi (y) - \sigma_\psi (y' ) | &\leq 4\lambda n^4 r + 
| \sigma_\varphi (x) - \sigma_\varphi (x' ) | \\
&\leq 4\lambda n^4 r + 2\| x-x' \| .
\end{align*}
Hence $\| y-y' \| \leq 8\lambda n^4 r + 4 \| x-x' \|$,
from which (iv) follows.

Finally, to prove (v) we suppose $x\geq 0$. By part (i) there is a
$y=(y_{ij})\in\mathcal{N}^\lambda_{L^n ,Y}(x)$. We then have, for $1\leq i,j
\leq n$,
\begin{align*}
L(x_{ij}, (\re (y))_{ij}) &\leq \frac12 L^e (x_{ij},y_{ij}) + \frac12 L^e 
(x^*_{ji}, y^*_{ji}) \\
&= \frac12 L^e (x_{ij},y_{ij}) + \frac12 L^e (x_{ji}, y_{ji}) \\
&\leq L^n (x,y)
\end{align*}
using the adjoint invariance of $L^e$, and so $\re (y)$ is a self-adjoint
element of $\mathcal{N}^\lambda_{L^n ,Y}(x)$. Suppose now that $y$ is
an arbitrary self-adjoint element of $\mathcal{N}^\lambda_{L^n ,Y}(x)$.
If $\psi\in SCP_n (Y)$ then as in the proof of (ii) there is a $\varphi\in
SCP_n (X)$ such that $| \sigma_\varphi (x) - \sigma_\psi (y) | \leq 
2\lambda n^4 r$, and thus since $\sigma_\varphi (x) \geq 0$ and $y$ is
self-adjoint we infer that
$\sigma_\varphi (y) \geq -2\lambda n^4 r$. Hence we conclude that $y\geq
-2\lambda n^4 r$.
\end{proof}


\begin{proposition}\label{P-orderisom}
Let $(X,L_X )$ and $(Y,L_Y )$ be Lip-normed operator systems, and 
suppose $\dist^n_s (X,Y) = 0$ for some $n\in\mathbb{N}$. 
Then there is a unital order isomorphism $\Phi : M_n \otimes X \to M_n 
\otimes Y$, and in the case $n=1$ we may arrange that 
$\mathcal{D}(L_Y ) = \Phi (\mathcal{D}(L_X ))$ and
$L_Y (\Phi (x)) = L_X (x)$ for all $x\in\mathcal{D}(L_X )$.
\end{proposition}

\begin{proof}
By assumption there is a sequence $\{ L_k \}_{k\in\mathbb{N}}$ of 
Lip-norms in $\mathcal{M}(L_A , L_B )$ such that $\lim_{k\to\infty}r_k = 0$ 
where $r_k = \dist^{\rho_{L_k ,n}}_H (UCP_n (X) , UCP_n (Y))$.
Let $x\in\mathcal{D}(L^n_X )$ and 
$\lambda > 2 L^n_X (x)$. Set $s = 4(\| x \| + 2\lambda n^4 )$.
In view of Lemma~\ref{L-diambound}(ii) we may assume (by removing
finitely many of the $L_k$'s and reindexing the sequence if necessary) 
that the sets $\mathcal{N}^\lambda_{L^n_k ,X}(x)$ for $k\geq 1$ are all 
contained in $\mathcal{D}_s (L^n_Y ) \cap\mathcal{B}^{M_n \otimes Y}_s$. 
This latter set is compact, since it is closed by the lower semicontinuity
of $L^n_Y$ and for any $t>0$ the set 
$\mathcal{D}_t (L_Y ) \cap\mathcal{B}^Y_t$ is 
compact (see the second half of the proof of Proposition~\ref{P-ntopology}), 
which implies the compactness of the set 
$\mathcal{D}_t (L^e_Y ) \cap\mathcal{B}^Y_t$ 
(use the decomposition of elements into real and self-adjoint parts) and
hence also the total boundedness of the tensor product of
$\mathcal{D}_t (L^e_Y ) \cap\mathcal{B}^Y_t$ with $M_n$. 
Since by Lemma~\ref{L-diambound}(iii) the 
diameters of $\mathcal{N}^\lambda_{L^n_k ,X}(x)$ converge to zero as $k\to
\infty$, this implies the existence of a subsequence of
$\big\{ \mathcal{N}^\lambda_{L^n_k ,X}(x) \big\}_{k\in\mathbb{N}}$ which 
converges in Hausdorff distance to some singleton $\{ \Phi (x) \}$. This
singleton must in fact be the same for each $\lambda > 2 L^n_X (x)$ because
for each $k\in\mathbb{N}$ we have
$$ \mathcal{N}^\lambda_{L^n_k ,X}(x) \subset 
\mathcal{N}^{\lambda'}_{L^n_k ,X}(x) $$
whenever $\lambda\leq\lambda'$.
Using a diagonal argument and relabeling indices we may assume that,
for all $x$ in a countable dense subset $D$ of $\mathcal{D}(L^n_X )$,
the sets $\mathcal{N}^\lambda_{L^n_k ,X}(x)$ for $\lambda > 2L^n_X (x)$
converge in Hausdorff distance as $k\to\infty$ to some singleton 
$\{ \Phi (x) \}$. Then in fact for any $x\in\mathcal{D}(L^n_X )$ and
$\lambda > 2L^n_X (x)$ the sets $\mathcal{N}^\lambda_{L^n_k ,X}(x)$
converge as $k\to\infty$ to some 
singleton $\{ \Phi (x) \}$, since for any $\epsilon > 0$ we can take
an $x' \in D$ with $\| x - x' \| \leq\epsilon /16$ and $2L^n_X (x' ) <
\lambda$ (since we may assume that $D$ was chosen so that $D\cap
\mathcal{D}_q (L^n_X )$ is dense in $\mathcal{D}_q (L^n_X )$ for all
positive rational numbers $q$) and a $k_0 \in\mathbb{N}$
such that, for all $k\geq k_0$, $\mathcal{N}^\lambda_{L^n_k ,X}(x' )$
is within Hausdorff distance $\epsilon /2$ of $\{ \Phi (x' ) \}$ and 
$8\lambda n^4 r_k \leq\epsilon /2$, from which it can be seen using 
Lemma~\ref{L-diambound}(iv) that 
for all $k\geq k_0$ the set $\mathcal{N}^\lambda_{L^n_k ,X}(x)$
lies inside the ball of radius $\epsilon$ centred at $\Phi (x' )$.

Now if $\mu\in\mathbb{C}$ and $x,x' \in\mathcal{D}(L^n_X )$ then for 
$\lambda > 2\big( |\mu | L^n_X (x) + L^n_X (x')\big)$ it is 
easily seen that
$$ \mathcal{N}^\lambda_{L^n_k , Y}(\mu x + x' ) \supset
\{ \mu y + y' : y\in\mathcal{N}^\lambda_{L^n_k , Y}(x)\text{ and }
y' \in\mathcal{N}^\lambda_{L^n_k , Y}(x' ) \} $$
so that $\mathcal{N}^\lambda_{L^n_k , Y}(\mu x + x' )$ converges to
$\{ \mu\Phi (x) + \Phi (x' ) \}$ as $k\to\infty$. Also, for every
$x\in\mathcal{D}(L^n_X )$ and $\lambda > L^n_X (x)$ we have 
$y^* \in\mathcal{N}^\lambda_{L^n_k , Y}(x^* )$ if and only if
$y\in\mathcal{N}^\lambda_{L^n_k , Y}(x)$, so that $\{ \Phi (x^* ) \}$ is
equal to the limit of $\{ y^* : y\in\mathcal{N}^\lambda_{L^n_k , Y}(x)\}$,
which by the continuity of the involution must be $\{ (\Phi (x))^* \}$. 
Hence we have
defined a $^*$-linear map $\Phi : \mathcal{D}(L^n_X ) \to M_n \otimes Y$. 
Note also that $\Phi$ is unital since $1_Y \in
\mathcal{N}^\lambda_{L^n_k , Y}(1_X )$ for all $\lambda > 0$.
We furthermore have
by Lemma~\ref{L-diambound}(ii) that the norm of $\Phi$ is bounded by
$4$ on $D$, and thus, since $D$ is dense in $\mathcal{D}(L^n_X )$ 
which in turn is dense in $M_n \otimes X$ by Lemma~\ref{L-dense}, 
$\Phi$ extends uniquely to a bounded $^*$-linear map from $M_n \otimes X$ 
to $M_n \otimes Y$, which we will again denote by $\Phi$.

Now by another diagonal argument and index relabeling we 
may assume that $\mathcal{N}^\lambda_{L^n_k , X}(y)$ converges in
Hausdorff distance as $k\to\infty$ to a singleton $\{ \Gamma (y) \}$ for 
all $y$ in a countable dense subset of $\mathcal{D}(L^n_Y )$ which contains 
$\Phi (D)$. We thus obtain, as above, a bounded unital $^*$-linear map
$\Gamma : \mathcal{D}(L^n_Y ) \to\mathcal{D}(L^n_X )$. We will show that
$\Phi$ and $\Gamma$ are mutual inverses. Suppose 
$x\in\mathcal{D}(L)$ and $\lambda > 2\max (L^n_X (x) , L^n_Y (\Phi (x)))$. 
For each $k\in
\mathbb{N}$ choose $x_k' \in\mathcal{N}^\lambda_{L^n_k , X}(\Phi (x))$
and $y_k \in\mathcal{N}^\lambda_{L^n_k , Y}(x)$. Let $\epsilon > 0$ and
$\varphi\in SCP_n (Y)$. Pick $k_0 \in\mathbb{N}$ large enough so that, for
all $k\geq k_0$, $2\lambda n^4 r_k \leq\epsilon$
and $\| y_k - \Phi (x) \| \leq\epsilon$. Then as in the proof of
Lemma~\ref{L-diambound}(i) for any $k\geq k_0$ we can 
find a $\psi\in SCP_n (X)$ such that
$$ | \sigma_\varphi (x) - \sigma_\psi (y_k ) | \leq 2n^4 L^e (x,y_k )r_k
\leq 2n^4 \lambda r_k \leq\epsilon $$
and similarly $| \sigma_\varphi (x_k' ) - \sigma_\psi (\Phi (x)) | \leq
\epsilon$, whence by the triangle inequality
\begin{align*}
| \sigma_\varphi (x) - \sigma_\varphi (x_k' ) | &\leq
2\epsilon + | \sigma_\psi (y_k ) - \sigma_\psi (\Phi (x)) | \\
&\leq 3\epsilon .
\end{align*}
Therefore $\| x - x_k' \| \leq 6\epsilon$, and so we have
$\lim_{k\to\infty}x_k' = x$. Hence
$$ \Gamma (\Phi (x)) = \Gamma \Big( \lim_{k\to\infty}y_k \Big) =
\lim_{k\to\infty}\Gamma (y_k ) = \lim_{k\to\infty} x_k' = x . $$
By a similar argument $\Phi (\Gamma (y)) = y$ for all $y\in
\mathcal{D}(L^n_Y )$, and hence by continuity we conclude 
that $\Phi$ and $\Gamma$ are mutual inverses.

Next we show that $\Phi$ and $\Gamma$ are positive. If $x\in\mathcal{D}
(L^n_X )$, $x\geq 0$, and $\lambda > 2L^n_X (x)$, 
then Lemma~\ref{L-diambound}(v) yields, for all $k\in\mathbb{N}$, a $y_k 
\in\mathcal{N}^\lambda_{L^n_k , Y}(x)$ with $y_k \geq -2\lambda n^4 r_k$.
Then $y$ is the limit as $k\to\infty$ of the positive elements 
$y_k + 2\lambda n^4 r_k$ and hence $y$ itself is positive. Thus $\Phi$ is
positive, and by a symmetric argument so is $\Gamma$. Hence $\Phi$ is a
unital order isomorphism.

It remains to show that if $n=1$ then $\Phi$ is isometric with respect to 
$L_X$ and $L_Y$, that is, $\mathcal{D}(L_Y ) = \Phi (\mathcal{D}(L_X ))$ and
$L_Y (\Phi (x)) = L_X (x)$ for all $x\in\mathcal{D}(L_X )$. 
Let $x\in\mathcal{D}(L_X )$ and 
$\lambda > \max (1,2 L_X (x))$. Suppose $\sigma , \sigma' \in
S(Y)$, and let $\epsilon > 0$. Choose $k\in\mathbb{N}$ large enough
so that we can find $\omega , \omega' \in S(X)$ with
$\lambda\rho_{L_k ,1}(\sigma , \omega ) \leq\epsilon$ and 
$\lambda\rho_{L_k ,1}(\sigma' , \omega' ) \leq\epsilon$, as well as a
$y_k \in\mathcal{N}^\lambda_{L^1_k ,X}(x)$ with 
$\| \Phi (x) - y_k \| \leq\epsilon$. Then
\begin{align*}
| \sigma (\Phi (x)) - \omega (x) | &\leq | \sigma 
(\Phi (x)) - \sigma (y_k ) | + | \sigma (y_k ) - \omega (x) | \\
&\leq \epsilon + \lambda \rho_{L_k ,1}(\sigma , \omega ) \\
&\leq 2\epsilon 
\end{align*}
and similarly $| \sigma' (\Phi (x)) - \omega' (x) | \leq
2\epsilon$. Thus, since 
\begin{align*}
\rho_{L_Y ,1}(\omega , \omega' ) &\leq 
\rho_{L_k ,1}(\omega , \sigma ) + \rho_{L_X ,1}(\sigma , \sigma' ) +
\rho_{L_k ,1}(\sigma' , \omega' ) \\ 
&\leq \rho_{L_Y ,1}(\sigma , \sigma' ) + 2\epsilon ,
\end{align*}
we have
\begin{align*}
| \sigma (\Phi (x)) - \sigma' (\Phi (x)) | &\leq
| \sigma (\Phi (x)) - \omega (x) | + 
| \omega (x) - \omega' (x) | +
| \omega' (x) - \sigma' (\Phi (x)) | \\
&\leq 4\epsilon + \rho_{L_X , 1}(\omega , \omega' ) L_X (x) \\
&\leq 2\epsilon (2 + L_X (x)) + \rho_{L_Y ,1}(\sigma , \sigma' )
L_X (x) .
\end{align*}
Dividing by $ \rho_{L_Y , 1}(\sigma , \sigma' )$ and letting $\epsilon\to
0^+$, we conclude that $L_Y (\Phi (x)) \leq L_X (x)$.
Since the above argument also applies to $\Gamma$ we must in fact have
$\mathcal{D}(L_Y ) = \Phi (\mathcal{D}(L_X ))$ and
$L_Y (\Phi (x)) = L_X (x)$ for all $x\in\mathcal{D}(L_X )$. 
\end{proof}

We will also need to know, for the proof of Theorem~\ref{T-orderisom}, 
that zero $m$-distance implies zero $n$-distance
for $m>n\geq 1$, which is a consequence of the following lemma.
 
\begin{lemma}\label{L-mnbound}
If $(X,L_X )$ and $(Y,L_Y )$ are Lip-normed operator systems and 
$m > n\geq 1$, then $\dist^n_s (X,Y) \leq\dist^m_s (X,Y)$.
\end{lemma}

\begin{proof}
Let $L\in\mathcal{M}(L_X , L_Y )$ and $\varphi\in UCP_n (X)$. Let
$\omega$ be any state on $X$, and define $\tilde{\varphi}\in 
UCP_m (X)$ by setting $\varphi' (x) = \varphi (x) + \omega (x)p$ 
for all $x\in X$, where $M_n$ has been identified with the
upper left-hand corner of $M_m$ and $p$ is the unit for
the lower right $(m-n) \times (m-n)$ corner. Choose $\psi' \in
UCP_m (Y)$ with $\rho_{L,m} (\varphi' , \psi' ) \leq
\dist_s^m (X,Y)$. If $\psi$ is the cut-down of $\varphi'$
to the upper-left hand $n\times n$ corner of $M_m$, then viewing it as an 
element of $UCP_n (X)$ we evidently have 
$$ \rho_{L,n} (\varphi , \psi ) \leq\rho_{L,m} (\varphi' , \psi' ) . $$
Hence $\dist_H^{\rho_{L,n}} (UCP_n (X) , UCP_n (Y))\leq
\dist_H^{\rho_{L,m}} (UCP_m (X) , UCP_m (Y))$, and so we conclude that
$\dist^n_s (X,Y) \leq\dist^m_s (X,Y)$.
\end{proof}

\begin{theorem}\label{T-orderisom}
Let $(X,L_X )$ and $(Y,L_Y )$ be Lip-normed operator systems.
\begin{enumerate}
\item[(i)] If $n\in\mathbb{N}$ then $\dist^n_s (X,Y) = 0$
if and only if there is a bi-Lip-isometric unital $n$-order isomorphism
between $X$ and $Y$.
\item[(ii)] We have $\dist_s (X,Y) = 0$ if and only if there is a 
bi-Lip-isometric unital complete order isomorphism between $X$ and $Y$.
\end{enumerate}
\end{theorem}

\begin{proof}
Proposition~\ref{P-Lipisometric} takes care of the ``if'' in each
part, and so we need only worry about the ``only if'' direction. 
Suppose $\dist^n_s (X,Y) = 0$ for some $n\in\mathbb{N}$.
Then by Lemma~\ref{L-mnbound} we have $\dist^m_s (X,Y) = 0$
for each $m=1, \dots ,n$. Applying the proof of 
Proposition~\ref{P-orderisom} successively for each $m=1, \dots ,n$ so 
that we can apply a diagonal argument across these values of $m$, we 
can find, for each $m=1, \dots ,n$, a unital order isomorphism
$\Phi_m : M_m \otimes X \to M_m \otimes Y$ such that, for each
$x\in M_m \otimes X$ and $\lambda > 2 L^m_X (x)$, 
the singleton $\{ \Phi (x) \}$ is the limit
of $\mathcal{N}^\lambda_{L^m_k , Y}(x)$ as $k\to\infty$ with respect to 
Hausdorff distance for a sequence of Lip-norms 
$L_k \in\mathcal{M}(L_A , L_B )$. By Proposition~\ref{P-orderisom} we
may also arrange that $\mathcal{D}(L_Y ) = \Phi (\mathcal{D}(L_X ))$ 
and $L_Y (\Phi_1 (x)) = L_X (x)$ for all $x\in\mathcal{D}(L_X )$.
We will show that $\Phi_m = \id_m \otimes\Phi_1$ for each $m=2,\dots ,n$.
Suppose then that $x\in\mathcal{D}(L_X )$. For each $k\in\mathbb{N}$
choose $y_k \in\mathcal{N}^\lambda_{L^1_k , Y}(x)$. If $e_{ij}$ is a
standard matrix unit in $M_m$ then 
$L^m_k (e_{ij} \otimes x , e_{ij} \otimes y_k ) = L^e_k (x,y_k )$
so that $e_{ij} \otimes y_k \in\mathcal{N}^\lambda_{L^m_k , Y}
(e_{ij} \otimes x)$. Thus 
$\mathcal{N}^\lambda_{L^m_k , Y}(e_{ij} \otimes x)$ must converge 
to the singleton containing
$$ \lim_{k\to\infty} e_{ij} \otimes y_k = e_{ij} \otimes\lim_{k\to\infty} 
y_k = (\id_m \otimes\Phi_1 )(e_{ij} \otimes x) , $$
whence $\Phi_m (e_{ij} \otimes x) = (\id_m \otimes\Phi_1 )(e_{ij} \otimes 
x)$. Since by Lemma~\ref{L-dense} the span of elements of the form 
$e_{ij} \otimes x$ with
$x\in\mathcal{D}(L^e_X )$ is dense in $M_m \otimes X$, we conclude that
$\Phi_m = \id_m \otimes\Phi_1$, so that $\Phi_1$ is a bi-Lip-isometric
$n$-order isomorphism. We thus obtain (i).

For (ii) we can use essentially the same proof, with the diagonal
arguments now extended across all $n\in\mathbb{N}$.
\end{proof}
 
\begin{corollary}
Let $A$ and $B$ be unital $C^*$-algebras with Lip-norms $L_A$ and
$L_B$, respectively.
\begin{enumerate}
\item[(i)] We have $\dist^1_s (A,B) = 0$ if and only if
there is a bi-Lip-isometric unital order isomorphism between $A$ and $B$.
\item[(ii)] If $n\geq 2$ then $\dist^n_s (A,B) = 0$
if and only if there is a bi-Lip-isometric $^*$-isomorphism between
$A$ and $B$.
\item[(iii)] We have $\dist_s (A,B) = 0$ if and only if there is a 
bi-Lip-isometric $^*$-isomorphism between $A$ and $B$.
\end{enumerate}
\end{corollary}

\begin{proof}
The corollary is an immediate consequence of Theorem~\ref{T-orderisom}
and the fact that a unital $2$-order isomorphism between $A$ and $B$ is
automatically a $^*$-isomorphism \cite{Ch}.
\end{proof} 

We remark that a unital order isomorphism between unital $C^*$-algebras 
need not be a $^*$-isomorphism. For instance, a unital $C^*$-algebra $A$ is 
always unitally order isomorphic to its opposite algebra $A^{\op}$, but 
these need not be $^*$-isomorphic, as the examples in \cite{Ph} 
demonstrate.

\section{$f$-Leibniz complete distance and convergence}\label{S-complete}

Let $(\mathcal{R} , \dist_s )$ be the metric space, under complete distance, 
of equivalence classes of Lip-normed operator systems with respect to 
bi-Lip-isometric unital complete order isomorphism. For economy
will simply refer to the elements of $\mathcal{R}$ as Lip-normed operator 
systems. Let $f : \mathbb{R}^4_+ \to\mathbb{R}_+$ be a continuous function. 
Given a Lip-normed unital $C^*$-algebra $(A,L)$, the 
Lip-norm $L$ is said to be {\it $f$-Leibniz} if it is the restriction of an
adjoint-invariant semi-norm $L'$ on $A$ which is finite on a dense 
$^*$-subalgebra and satisfies the {\it $f$-Leibniz property}
$$ L'(xy) \leq f(L'(x), L'(y), \| y \| , \| x \| ) $$
for all $x,y\in A$ with $L' (x), L' (y) < \infty$. 
When $f(a,b,c,d) = ac + bd$ this is the usual Leibniz rule, in which case 
we simply say that $L$ is {\it Leibniz}.
The Lip-normed unital $C^*$-algebras of Example~\ref{E-actions} are 
Leibniz, as are those obtained from Lipschitz semi-norms on functions 
over a compact metric space.
We denote by $\mathcal{R}_{\alg}$ the subset of $\mathcal{R}$ consisting
of Lip-normed unital $C^*$-algebras, and for $(A, L_A )$ and
$(B, L_B )$ in $\mathcal{R}_{\alg}$ we define the {\it $f$-Leibniz complete 
distance} $\dist_{s,f}(A,B)$ in the same way that the complete distance is 
defined (Definition~\ref{D-dist}) except that the infimum is now taken
over the $f$-Leibniz Lip-norms in $\mathcal{M}(L_A , L_B )$ (if no such 
$f$-Leibniz Lip-norm exists we set $\dist_{s,f}(A,B) = \infty$). Note that
$\dist_{s,f}$ might not satisfy the triangle inequality without further 
hypotheses on $f$, but that will not be of consequence for our application 
here, and we can still speak of Cauchy sequences with respect to 
$\dist_{s,f}$ in the obvious sense. It can be
seen that the estimates in Example~\ref{E-sphere} for complete distance
also apply to $f$-Leibniz complete distance for suitable $f$ (although
$f$ may depend on the matrix algebra), and if $N$ is a bridge between
two Leibniz Lip-normed $C^*$-algebras of the form that appears in 
Proposition~\ref{P-distest} then the resulting Lip-norm on the direct sum 
is Leibniz (see Section~\ref{S-tb} for examples of the use of 
bridges like those in Proposition~\ref{P-distest}). In this section we show 
that every sequence in $\mathcal{R}_{\alg}$ which is
Cauchy with respect to $f$-Leibniz complete distance converges in 
$\mathcal{R}_{\alg}$ with respect to complete distance. We may
interpret this as saying that the metric space $(\mathcal{R}_{\alg},
\dist_s )$ is ``complete'' relative to $f$-Leibniz complete distance.

I would like to thank Narutaka Ozawa for suggesting the idea behind 
the proof of the following lemma. 
Given a sequence $\{ A_k \}_{k\in\mathbb{N}}$ of $C^*$-algebras 
we denote by $\prod A_k$ the $C^*$-algebra of bounded sequences with the 
supremum norm and by $\bigoplus A_k$ the $C^*$-subalgebra of sequences 
converging to zero. 

\begin{lemma}\label{L-statelim}
Let $\{ A_k \}_{k\in\mathbb{N}}$ be a sequence of unital $C^*$-algebras
and $X$ a separable operator subsystem of $\prod A_k \big/ \bigoplus A_k$,
and let $n\in\mathbb{N}$. Suppose that, for each $x\in M_n \otimes X$, at 
least one lift (and hence every lift) $\sum e_{ij} \otimes (x_k^{ij} )_k$  
of $x$ with respect to the quotient $\big( M_n \otimes\big( \prod A_k \big) 
\big) \big/ \big( M_n \otimes \big( \bigoplus A_k \big) \big) \cong
M_n \otimes\big( \prod A_k \big/ \bigoplus A_k \big)$ satisfies
$$ \lim_{k\to\infty} \Big\| \sum e_{ij} \otimes x_k^{ij} \Big\| = \| x \| . $$
Then for every $\varphi\in UCP_n (X)$ there are
$\varphi_k \in UCP_n (A_k )$ for $k\in\mathbb{N}$ such that for all 
$(x_k )_k + \bigoplus A_k \in X$ we have
$$ \varphi \big( (x_k )_k + {\textstyle\bigoplus} A_k \big) = 
\lim_{k\to\infty}\varphi_k (x_k ) . $$
\end{lemma}

\begin{proof}
First we consider an arbitrary finite-dimensional operator subsystem 
$Y$ of $X$ and show that the conclusion of the lemma holds with respect to
elements of $Y$.
Letting $\pi : \prod A_k \to\prod A_k \big/ \bigoplus A_k$ be the quotient 
map, there exists, by elementary linear algebra, a unital linear map 
$x\stackrel{\alpha}{\mapsto}\alpha (x) = (\alpha (x)_k )_k$ from $X$ to 
$\prod A_k $ such that $\pi\circ\alpha = \id_X$. We may
assume that $\alpha$ is Hermitian for otherwise we can replace it with 
its real part $(\alpha + \alpha^* )/2$. Since $Y$ is finite-dimensional
the unit ball of $M_n \otimes Y$ is compact, and so by our assumption
on lifts of elements of $M_n \otimes X$ 
we can find a sequence $\delta_1 , \delta_2 , \dots$ of positive real
numbers with $\lim_{k\to\infty} \delta_k = 0$
such that, for all $\sum e_{ij} \otimes x_{ij}$ in  
the unit ball of $M_n \otimes Y$ and $k\in\mathbb{N}$,
$$ \Big\| \sum e_{ij} \otimes x_{ij} \Big\| - \delta_k < 
\Big\| \sum e_{ij} \otimes \alpha (x_{ij})_k \Big\| < \Big\| 
\sum e_{ij} \otimes x_{ij} \Big\| + \delta_k . $$
This implies in particular that for each sufficiently large 
$k\in\mathbb{N}$ the map $\pi_k \circ\alpha$ is injective on $Y$, 
where $\pi_k : \prod A_k \to A_k$ is the 
projection map. Let $\varphi\in UCP_n (Y)$. For each sufficiently large
$k\in\mathbb{N}$ we can define the linear map $\psi_k : (\pi_k \circ\alpha )
(Y) \to M_n$ by
$$ \psi_k (a) = \varphi ((\pi_k \circ\alpha )^{-1}(a)) $$
for all $a\in (\pi_k \circ\alpha )(Y)$.
Then $\psi_k$ is unital and Hermitian, and $\| \id_n \otimes\psi_k \|
\leq (1-\delta_k )^{-1}$. 
By \cite[Thm.\ 2.10]{Sm} the completely bounded norm
$\| \psi_k \|_{cb}$ is equal to $\| \id_n \otimes\psi_k \|$ and hence is at 
most $(1-\delta_k )^{-1}$. By the Wittstock extension theorem 
(see \cite{Was})
there is an extension of $\psi_k$ to $A_k$ with the same completely
bounded norm. We denote this extension also by $\psi_k$.
By the Wittstock decomposition theorem (see \cite{Was}) there exist
completely positive maps $\psi^+_k$ and $\psi^-_k$ from $A_k$ to
$M_n$ such that $\psi_k = \psi^+_k - \psi^-_k$ and 
$\| \psi_k \|_{cb} \geq \| \psi^+_k + \psi^-_k \|$. We then have
$$ \psi^+_k (1) = \psi_k (1) + \psi^-_k (1) = 1 + \psi^-_k (1) \geq 1 $$
and
$$ \| \psi^+_k (1) \| = \| \psi^+_k \| \leq 
\| \psi^+_k + \psi^-_k \| \leq \| \psi_k \|_{cb} \leq (1-\delta_k )^{-1} . $$
Also,
$$ \| \psi^+_k - \psi_k \| = \| \psi^+_k (1) - \psi_k (1) \| 
= \| \psi^+_k (1) - 1 \| \leq (1-\delta_k )^{-1} - 1 . $$
Since $\psi^+_k (1) > 0$ we can define the u.c.p.\ map $\varphi_k :
A_k \to M_n$ by 
$$ \varphi_k (a) = \psi^+_k (1)^{-\frac12} \psi^+_k (a) 
\psi^+_k (1)^{-\frac12} $$
for all $a\in A_k$. Then
\begin{align*}
\| \psi_k - \varphi_k \| &\leq \| \psi_k - \psi^+_k \| +
\| \psi^+_k - \varphi_k \| \\
&\leq ((1-\delta_k )^{-1} - 1) + \| 1-\psi^+_k (1)^{-\frac12} \| 
\| \psi^+_k \| (1 + \| \psi^+_k (1)^{-\frac12} \|) \\
&\leq ((1-\delta_k )^{-1} - 1) + 2(1-(1-\delta_k )^{\frac12})
(1-\delta_k )^{-1} , 
\end{align*} 
and this last expression tends to zero as $k\to\infty$.
It follows that, for all $(x_k )_k + \bigoplus A_k \in Y$,
$$ \varphi \big( (x_k )_k + {\textstyle\bigoplus} A_k \big) = 
\lim_{k\to\infty}\psi_k (x_k ) = \lim_{k\to\infty}\varphi_k (x_k ) . $$

Now suppose that $X_1 \subset X_2 \subset\cdots$ is an
increasing sequence of finite-dimensional operator subsystems of $X$
with union dense in $X$. Let $\varphi\in UCP_n (X)$. 
Then for each $j\in\mathbb{N}$ there exists by the first paragraph 
u.c.p.\ maps $\varphi_k$ on $A_k$ for sufficiently large $k$ (and hence
for all $k$) such that 
$$ \varphi \big( (x_k )_k + {\textstyle\bigoplus} A_k \big) = 
\lim_{k\to\infty}\varphi_k (x_k ) $$
for all $(x_k )_k + \bigoplus A_k \in X_j$. By applying a diagonal
argument over $j\in\mathbb{N}$ we can assume that the equality in the
above display holds for all $(x_k )_k + \bigoplus A_k \in
\bigcup_{j\in\mathbb{N}}X_j$. A straightforward approximation argument then
shows that this equality in fact holds for all $(x_k )_k + \bigoplus A_k
\in X$, completing the proof.
\end{proof}

\begin{lemma}\label{L-quotientest}
Let $(Z,L_Z )$ be a Lip-normed operator system, $\Phi : Z\to X$ and
$\Gamma : Z\to Y$ surjective u.c.p.\ maps onto operator systems $X$ and
$Y$, and $L_X$ and $L_Y$ the quotient Lip-norms induced via $\Phi$ and
$\Gamma$, respectively. Then
$$ \dist_s (X,Y)\leq\sup_{n\in\mathbb{N}}\,\dist_H^{\rho_{L,n}} 
(UCP_n (X), UCP_n (Y)) , $$
with $UCP_n (X)$ and $UCP_n (Y)$ considered as subsets of $UCP_n (Z)$.
\end{lemma}

\begin{proof}
Set $r=\sup_{n\in\mathbb{N}}\dist_H^{\rho_{L,n}} (UCP_n (X), UCP_n (Y))$
(as can be seen from the proof of Proposition~\ref{P-completediambound},
this supremum is bounded by $\diam (Z,L_Z )$).
As in \cite[Example 5.6]{GHDQMS} for any $\gamma > 0$ we can construct a 
bridge between two copies of $(Z,L_Z )$ by setting $N(z,z' ) = \gamma^{-1}
\| z - z' \|$. Let $M$ be the Lip-norm 
$$ L(z,z' ) = \max (L_Z (z) , L_Z (z') , N(z,z' )) $$
on $\mathcal{D}(L_Z ) \oplus\mathcal{D}(L_Z )$ and $L$ the quotient
Lip-norm induced by $M$ via the u.c.p.\ map $(z,z' )\mapsto (\Phi (z),
\Gamma (z' ))$. Then $L\in\mathcal{M}(L_X , L_Y )$. Denote the 
projections of $Z\oplus Z$ onto the first and second direct summand by
$\pi_1$ and $\pi_2$, respectively. 

Now suppose $\varphi\in UCP_n (X)$. Then by assumption for some 
$\psi\in UCP_n (Y)$ we have $\rho_{L_Z ,n} (\varphi\circ\Phi , \psi\circ
\Gamma ) \leq r$. Also, if $(z,z' )\in\mathcal{D}_1 (M)$ then $\| z-z' \| 
\leq\gamma$ so that $\| (\psi\circ\Gamma )(z) - (\psi\circ\Gamma )
(z' ) \| \leq\gamma$ and hence $\rho_{M,n} 
(\psi\circ\Gamma\circ\pi_1 , \psi\circ\Gamma\circ\pi_2 )\leq\gamma$
(where to avoid confusion we have included the compositions with
the projection maps $\pi_1$ and $\pi_2$, contrary to our usual
practice). Thus, since $\rho_{L,n}$ is the restriction of 
$\rho_{M,n}$ via the identification arising from the quotient map, we have 
by the triangle inequality
\begin{align*}
\rho_{L,n}(\varphi , \psi ) &\leq \rho_{M,n}(\varphi\circ\Phi , \psi\circ
\Gamma ) + \rho_{M,n}(\psi\circ\Gamma\circ\pi_1 , 
\psi\circ\Gamma\circ\pi_2 ) \\
&\leq r + \gamma .
\end{align*}
Hence $\dist_s (X,Y)\leq r + \gamma$, which yields the result since
$\gamma$ was arbitrary.
\end{proof}

In the proof of the following theorem, we will abbreviate expressions of the
form $\dist^{\rho_{L,n}}_H (UCP_n (X), UCP_n (Y))$ to 
$\rho_{L,n} (UCP_n (X), UCP_n (Y))$ to reduce
the number of subscripts, and whenever we have a quotient Lip-norm then 
will identify the state space of the quotient operator system with a 
subset of the state space of the original operator system under the 
induced isometry (Proposition~\ref{P-isometry}) as is our usual practice
in the case of projections onto direct summands.

\begin{theorem}\label{T-complete}
Let $\{ (A_k , L_k ) \}$ be a sequence in $\mathcal{R}_{\alg}$ which is 
Cauchy with respect to $f$-Leibniz complete distance for a given
continuous $f : \mathbb{R}^4_+ \to\mathbb{R}_+$. Then 
$\{ (A_k , L_k ) \}$ converges in $\mathcal{R}_{\alg}$ with 
respect to complete distance.
\end{theorem}

\begin{proof}
To show that $\{ (A_k , L_k ) \}$ converges it suffices to show the 
convergence of a subsequence, and so we may assume that 
$\dist_{s,f} (A_k , A_{k+1}) < 
2^{-k}$ for all $k\in\mathbb{N}$. Then there exist $f$-Leibniz Lip-norms 
$L_{k,k+1} \in\mathcal{M}(L_k , L_{k+1})$ with
$$ \rho_{L_{k,k+1} ,n}(UCP_n (A_k ), UCP_n (A_{k+1})) < 2^{-k} $$
for all $n,k\in\mathbb{N}$. 
Let $Z$ be the set of sequences $( x_k )_k$ with $x_k \in\mathcal{D}
(L_k )$ such that, for some $\lambda > L_1 (x_1 )$, $x_{k+1}\in
\mathcal{N}^\lambda_{L_{k,k+1},A_{k+1}}(x_k )$ for all $k\in\mathbb{N}$
(see Definition~\ref{D-N}).
We will show that $Z$ is a subset of the direct product $\prod A_k$.
Let $J_k$ be the semi-norm on $\prod_{j=1}^k \mathcal{D}(L_j )$ given by
$$ J_k ((x_j )_j ) = \sup_{1\leq j \leq k-1}L_{j, j+1}(x_j , x_{j+1}) . $$
By \cite[Lemma 12.2]{GHDQMS} $J_k$ is a Lip-norm.
Denote by $Q_k$ the quotient Lip-norm on $\mathcal{D}(L_1 )\oplus 
\mathcal{D}(L_k )$ induced by $J_k$ via the projection map. Then 
$$ \rho_{Q_k ,1}(S(A_1 ), S(A_k )) \leq 2^{-1} + 2^{-2} + 
\cdots + 2^{-k} < 1 . $$
Suppose $( x_k )_k \in Z$, and let $\lambda > L_1 (x_1 )$ be such that
$x_{k+1}\in\mathcal{N}^\lambda_{L_{k,k+1},X_{k+1}}(x_k )$ for all 
$k\in\mathbb{N}$. Then, for each $k\in\mathbb{N}$, $x_k$ is an 
element of $\mathcal{N}^\lambda_{Q_k , A_k}(x_1 )$ and hence 
by Lemma~\ref{L-diambound} has norm bounded by 
$\| x_1 \| + 2\lambda\rho_{Q_k ,1}(S(A_1 ), S(A_k )) \leq 
\| x_1 \| + 2\lambda$. 
Therefore $( x_k )_k$ is a bounded sequence and so belongs 
to $\prod A_k$, as we wished to show.

We define the semi-norm $L_Z$ on $Z$ by
$$ L_Z (( x_k )_k ) = \sup_{k\in\mathbb{N}}L_{k,k+1}(x_k , x_{k+1}) $$ 
(which is finite by the definition of $Z$). 
Theorem 12.9 of \cite{GHDQMS} then shows that $L_Z$ is a Lip-norm on $Z$.
It can be seen by the $f$-Leibniz property that the operator system $B$ 
generated by $Z$ in $\prod A_k$ is in fact a $C^*$-algebra.
Let $A$ be the $C^*$-subalgebra of $\prod A_k / \bigoplus A_k$ which
is the image of $B$ under the quotient map $\pi : \prod A_k 
\to\prod A_k \big/ \bigoplus A_k$, and let $L$ be the quotient Lip-norm on 
$A$ induced by $L_Z$. Then $(A,L)$ is a Lip-normed unital $C^*$-algebra.
Our goal now is to show that $\{ (A_k , L_k ) \}$ 
converges to $(A,L)$ with respect to complete distance. 
By Lemma~\ref{L-quotientest} it suffices to show that,
for all $n\in\mathbb{N}$, $UCP_n (A)$ coincides with the Hausdorff limit
$H_n \subset UCP_n (B)$ of $\{ UCP_n (A_k ) \}_{k\in\mathbb{N}}$,
which exists due to the
completeness, in the Hausdorff metric, of the set of closed subspaces of 
the compact set $UCP_n (B)$. Note that for each $k\in\mathbb{N}$ the image
of $B$ under the projection onto $A_k$ is surjective so that indeed
$UCP_n (A_k )\subset UCP_n (B)$, and that the convergence
of $\{ UCP_n (A_k ) \}_{k\in\mathbb{N}}$ to $H_n$ is uniform over $n$
because the Cauchy condition is uniform over $n$ by assumption.

If $\{ \varphi_k \}_{k\in\mathbb{N}}$ is a sequence such that 
$\varphi_k \in UCP_n (A_k )$ and $\{ \varphi_k \circ\pi_k 
\}_{k\in\mathbb{N}}$ is point-norm convergent (necessarily to an element
of $H_n$), then setting
$$ \varphi \big( (x_k )_k + {\textstyle\bigoplus} A_k \big) = 
\lim_{k\to\infty}\varphi_k (x_k ) $$
for $(x_k )_k + \bigoplus A_k \in X$ we obtain a map 
$\varphi : X\to M_n$. This map is u.c.p.\ in view of 
the identification $M_n \otimes\big( \prod A_k \big/ \bigoplus 
A_k \big) \cong \big( M_n \otimes\big( \prod A_k \big) \big) \big/ \big(
M_n \otimes \big( \bigoplus A_k \big) \big)$ and the fact that positive 
elements in quotients lift to positive elements. We thus see that 
$H_n \subset UCP_n (A)$. 

It remains to show that $H_n \supset UCP_n (A)$. With a view to applying
Lemma~\ref{L-statelim}, we will show
that every $x\in\ M_n \otimes A$ has a lift $(x_k )_k \in M_n \otimes
\prod A_k$ satisfying $\lim_{k\to\infty} \| x_k \| = \| x \|$. 
Notice first that if $( z_k )_k
\in M_n \otimes Z$ then for some $\lambda$ not depending on $j$ we have
$\big| \| z_j \| - \| z_{j+1} \| \big| \leq 
2^{-j+1} n^4 \lambda$ by Lemma~\ref{L-diambound}(ii) (since each $z_j$ is 
self-adjoint), so that
$\{ \| z_k \| \}_{k\in\mathbb{N}}$ is a Cauchy sequence and hence
$\| \pi (( z_k )_k ) \| = \lim_{k\to\infty} \| z_k \|$. Now
suppose $x\in M_n \otimes A$ and let $(x_k )_k$ be a lift of
$x$ to $M_n \otimes B$. Then $(x_k^* x_k )_k \in M_n \otimes B$, and
so there exists a $(y_k )_k \in M_n \otimes Z$ such that 
$\| x_k^* x_k - y_k \| < \varepsilon$ for all $k\in\mathbb{N}$, and from 
above we have $\| \pi (( y_k )_k ) \| = \lim_{k\to\infty} \| y_k \|$.
Let $\varepsilon > 0$, and
choose $k_0 \in\mathbb{N}$ such that, for all $j,k \geq k_0$,
$\big| \| y_j \| - \| y_k \| \big| < \varepsilon$. Then, for all 
$j,k \geq k_0$,
\begin{align*}
\big| \| x_j \|^2 - \| x_k \|^2 \big| &= \big| \| x_j^* x_j \| -
\| x_k^* x_k \| \big| \\ 
&\leq \big| \| x_j^* x_j \| -
\| y_j \| \big| + \big| \| y_j \| - \| y_k \| \big| + \big| \| y_k \| -
\| x_k^* x_k \| \big| \\
&< 3\varepsilon .
\end{align*}
It follows that $\{ \| x_k \|^2 \}_{k\in\mathbb{N}}$ is a Cauchy sequence
and hence converges. Thus $\lim_{k\to\infty} \| x_k \|$ exists, and it must
equal $\| x \|$. We can therefore apply Lemma~\ref{L-statelim}, so that
given $\varphi\in UCP_n (A)$ there exist $\varphi_k \in UCP_n (A_k )$ 
for $k\in\mathbb{N}$ such that for all $(x_k )_k + \bigoplus
A_k \in A$ we have
$$ \varphi \big( (x_k )_k + {\textstyle\bigoplus} A_k \big) = 
\lim_{k\to\infty}\varphi_k (x_k ) , $$
whence $H_n \supset UCP_n (A)$. Thus $H_n$ and $UCP_n (A)$ coincide,
completing the proof.
\end{proof}

Using the arguments of this section we might hope to show that the
metric space $(\mathcal{R}, \dist_s )$ is complete. However,
without the sharp control on the norms of non-self-adjoint elements
that the $f$-Leibniz property provides in Theorem~\ref{T-complete},
we would not be able to apply Lemma~\ref{L-statelim}.

\section{Total boundedness}\label{S-tb}

We will establish a version of 
Theorem 13.5 in \cite{GHDQMS} (``the quantum Gromov compactness theorem'')
for complete distance using approximation by 
Lip-normed operator subsystems of matrix algebras. As before
$(\mathcal{R} , \dist_s )$ is the metric space 
of equivalence classes of Lip-normed operator systems with respect to 
bi-Lip-isometric unital complete order isomorphism. 

\begin{notation} 
For a Lip-normed operator system $(X,L)$ and 
$\epsilon > 0$ we denote by $\Afn_L (\epsilon )$ the smallest
integer $k$ such that there is a Lip-normed operator system 
$(Y,L_Y )$ with $Y$ an operator subsystem of the matrix algebra $M_k$ and 
$\dist_s (X,Y) \leq\epsilon$. If no such integer $k$ exists we write
$\Afn_L (\epsilon ) = \infty$. We denote by $\mathcal{R}_{\fa}$
the subset of $\mathcal{R}$ consisting of Lip-normed operator systems 
$(X,L)$ for which $\Afn_L (\epsilon )$ is finite for all $\epsilon > 0$.
\end{notation}

We remark that every Lip-normed nuclear operator system and 
Lip-normed unital exact $C^*$-algebra is contained
in $\mathcal{R}_{\fa}$ by Propositions~\ref{P-finiterank} and
\ref{P-exact}, respectively. Note also that $\mathcal{R}_{\fa}$ is
a closed subset of $\mathcal{R}$ under the complete distance topology.

\begin{lemma}\label{L-singleopsys}
Let $X$ be a finite-dimensional operator system and $C\geq 0$. Then the 
set $\mathcal{C} = \{ (X,L)\in\mathcal{R}$ : $\diam (X,L)\leq C \}$ is
totally bounded.
\end{lemma}

\begin{proof} 
Proposition 13.13 and the proof of Proposition 13.14 in \cite{GHDQMS} show 
that, given $\epsilon > 0$, there is a finite subset
$\mathcal{S}\subset\mathcal{C}$ such that for all $(X,L)\in
\mathcal{C}$ there is a $(X,L' )\in\mathcal{C}$ and a bridge
$N$ between $(X,L)$ and $(X,L' )$ of the 
form $N(x,y) = \epsilon^{-1}\| x - y \|$. Let $M$ be the Lip-norm in
$\mathcal{M}(L,L' )$ given by
$$ M(x,y) = \max (L (x) , L' (y) , N(x,y)) . $$
Now if $\phi\in UCP_n (X)$ and
$(x,y)\in\mathcal{D}_1 (M)$ then $\| x-y \| \leq\epsilon$
so that $\| \phi (x) - \phi (y) \| \leq\epsilon$, and so by
Proposition 2.10 we have $\rho_{M,n} (\phi\circ\pi_1 , \phi\circ\pi_2 )
\leq\epsilon$, where $\pi_1$ and $\pi_2$ are the projections of
$X\oplus X$ onto the first and second direct summands, respectively.
Hence $\dist_s ((X,L),(X,L'))\leq\epsilon$, from which we
conclude that $\mathcal{C}$ is totally bounded.
\end{proof}

\begin{theorem}\label{T-ccompact}
Let $\mathcal{C}$ be a subset of $\mathcal{R}_{\fa}$. Then $\mathcal{C}$
is totally bounded if and only if
\begin{enumerate}
\item[(i)] there is an $M > 0$ such that the diameter of every
element of $\mathcal{C}$ is bounded by $M$, and
\item[(ii)] there is a function $F: (0,\infty ) \to (0,\infty )$ with 
$\Afn_L (\epsilon ) \leq F(\epsilon )$ for all $(X,L) \in\mathcal{C}$ .
\end{enumerate}
\end{theorem}

\begin{proof}
For the ``only if'' direction, suppose that $\mathcal{C}$ is a totally 
bounded subset of $\mathcal{R}_{\fa}$. If there did not exist an $M>0$
bounding the complete diameter of every element of $\mathcal{C}$, then
we could find a sequence $\{ (X_k , L_k ) \}_{k\in\mathbb{N}}$ such that 
$\diam (X_{k+1} , L_{k+1}) \geq\diam (X_k , L_k )  + 1$ for every 
$k\in\mathbb{N}$,
in which case $\dist_s (X_k , X_{k'}) \geq 1$ for all $k,k' \in\mathbb{N}$,
contradicting total boundedness. To verify condition (ii), we can find a 
finite $(\epsilon /2)$-dense subset $\mathcal{G}$ of $\mathcal{C}$ and set
$$ G(\epsilon ) = \max\{ \Afn_L (\epsilon /2) : (X,L)\in\mathcal{G} \} . $$
Then by the triangle inequality $\Afn_L (\epsilon )\leq G(\epsilon )$
for any $(X,L)\in\mathcal{C}$.

To prove the converse, suppose that conditions (i) and (ii) hold. By 
(ii) we see that it is sufficient to prove, for $k\geq j\geq 1$ and $M > 0$, 
the total boundedness of the collection of Lip-normed 
operator systems $(X,L)$ where $X$ is an operator subsystem 
of the matrix algebra $M_k$ with $X_{\sa}$ of real linear dimension
$j$ (in which case we will say that $X$ has {\em Hermitian dimension $j$})
and $\diam (X,L) \leq M$. 
Since the closed unit ball of the self-adjoint part of $M_k$ is 
compact, the set of closed unit balls of the self-adjoint parts of
operator subsystems of $M_k$ of Hermitian dimension $j$ is totally bounded 
in the Hausdorff metric. Also, by Lemma~\ref{L-singleopsys},
for every $M > 0$ and operator subsystem $X$ of $M_k$ of Hermitian
dimension $j$ the set of Lip-normed operator
systems $(X,L)$ with $\diam (X,L)\leq M$ is totally bounded.
Thus we need only show that, for every $\epsilon > 0$ and $M>0$, 
if $X$ and $Y$ are operator subsystems of $M_k$ of Hermitian dimension $j$
the closed unit balls of the self-adjoint parts of which are within 
Hausdorff distance $(4k)^{-1}\epsilon\min (M^{-1} , 1)$, 
and $L_X$ is a Lip-norm on $X$ with $\diam (X,L_X )\leq M$, 
then there is a Lip-norm $L_Y$ on $Y$ with 
$\dist_s ((Y,L_Y ),(X, L_X )) \leq \epsilon$ and
$\diam (Y,L_Y ) \leq M + \epsilon$. So let $X$ and $Y$ be such 
operator systems and $L_X$ such a Lip-norm on $X$ for given 
$\epsilon > 0$ and $M>0$. We may assume that $\epsilon < 1/2$.
Set $\delta = (4k)^{-1}\epsilon\min (M^{-1} , 1)$.
By \cite[Lemma 3.2.3]{BK} there is a 
(real linear) projection $P$ from $(M_k )_{\sa}$ onto $Y_{\sa}$ of norm 
$\leq k$. The restriction $Q$ of $P$ 
to $X_{\sa}$ is a bijection, for if $x\in\cap X_{\sa}$ with $\| x \| = 1$
then we can find a $y\in Y_{\sa}$ with $\| y - x \| < \delta\leq
\epsilon /k$ so that
$$ \| Q(x) \| \geq \| y \| - \| Q(x - y) \| \geq 1 - \epsilon
> \frac12 , $$
yielding injectivity and hence also bijectivity since $X_{\sa}$ and 
$Y_{\sa}$ are of equal finite dimension. The above display also shows that
the norm of $Q^{-1}$ is bounded by $2$. We next define a semi-norm $L_X$ 
on $Q^{-1}(Y)$ by $L_X (x) = L_Y (Q(x))$ (note that $\mathcal{D}(L_Y )$ is
equal to $Y_{\sa}$ by finite-dimensionality). 
Since the restriction of $Q$ to $X$ is bijective and
$L_Y$ is a Lip-norm we must have $L_X (x) = 0$ if and only if $x\in
\mathbb{R}1$. Thus $L_X$ is a Lip-norm in view of the finite-dimensionality
of $X_{\sa}$, and $(X,L_X )$ is a Lip-normed operator system since
$\mathcal{D}(L_X ) = X_{\sa}$ and $\mathcal{D}_1 (L_X )$ is closed in 
$X_{\sa}$ by the bijectivity and continuity, respectively, of $Q$.

On $\mathcal{D}(L_X )\oplus\mathcal{D}(L_Y )$ we define the semi-norm 
$N$ by $N(x,y) = \epsilon^{-1}\| x - y \|$. We will argue that $N$ is a 
bridge. For this it suffices to show that, for all $x\in\mathcal{D}(L_X )$, 
$$ N(x,Q(x))\leq L_Y (Q(x)) , $$
for then in condition (ii) of Definition~\ref{D-bridge} given 
$x\in\mathcal{D}(L_X )$ we can take $Q(x)$, and given $y\in\mathcal{D}(L_Y )$
we can take $Q^{-1}(y)$. So let $x\in\mathcal{D}(L_X )$.
Then we can find a $y\in Y_{\sa}$ such that $\| x-y \| \leq\delta \| x \|$,
so that
$$ \| x - Q(x) \| \leq \| x-y \| + \| Q(x-y) \| \leq (1+k)\| x-y \|
\leq 2k\delta \| x \| \leq 4k\delta \| Q(x) \| . $$ 
Applying this estimate with $x$ replaced by $x-\lambda 1$ where $\lambda$
is the infimum of the spectrum of $Q(x)$, we have
\begin{align*}
\| x - Q(x) \| = \| x - \lambda 1 -
Q(x - \lambda 1) \| &\leq 4k\delta \| Q(x - \lambda 1) \| \\
&= 4k\delta \| Q(x) - \lambda 1 \| \\
&\leq 4kM\delta L_Y (Q(x) - \lambda 1) \\
&= 4kM\delta L_Y (Q(x)) \\
&\leq \epsilon L_Y (Q(x)) 
\end{align*}
with Proposition~\ref{P-normLip-norm} yielding the second inequality in the 
string. We thus conclude that $N$ is a bridge. Let $L$ be the Lip-norm 
in $\mathcal{M}(L_X , L_Y )$ given by
$$ L(x,y) = \max (L_X (x) , L_Y (y) , N(x,y)) . $$
It remains to show that
$\dist^{\rho_{L,n}}_H (UCP_n (X), UCP_n (Y)) \leq \epsilon$ for
all $n\in\mathbb{N}$, for then $\dist_s (X,Y) \leq\epsilon$ and hence
also $\diam (X,L_X )\leq\diam (Y,L_Y ) + \epsilon\leq M + \epsilon$. 
Let $\varphi\in UCP_n (X)$. By Arveson's extension theorem 
we can extend $\varphi$ to a u.c.p.\ map $\varphi' : M_k \to M_n$. We then
have
\begin{align*}
\rho_{L,n}(\varphi , \varphi' |_Y ) &= \sup \{ | \varphi' (x) - 
\varphi' (y) | : (x,y)\in\mathcal{D}_1 (L) \} \\
&\leq \sup \{ | \varphi' (x) - \varphi' (y) | : (x,y)\in X\oplus Y
\text{ and }\| x -y \| \leq\epsilon \} \\
&\leq\epsilon .
\end{align*}
Similarly, if $\varphi\in UCP_n (Y)$ then extending it by Arveson's
theorem to a u.c.p.\ map $\varphi' : M_k \to M_n$ we have 
$\rho_{L,n}(\varphi , \varphi' |_Y )\leq\epsilon$. Thus
$\dist^{\rho_{L,n}}_H (UCP_n (X), UCP_n (Y)) \leq \epsilon$, as desired.
\end{proof}

An immediate consequence of Theorem~\ref{T-ccompact} is  
the separability of $\mathcal{R}_{\fa}$.

\begin{corollary}
The metric space $\mathcal{R}_{\fa}$ is separable.
\end{corollary}

\begin{question}
Given $n>1$ and $M>0$, is the set of all $n$-dimensional Lip-normed
operator systems of diameter at most $M$ totally bounded and/or separable?
\end{question}

We may think of $\log\Afn_L (\epsilon )$ as an analogue of 
Kolmogorov $\epsilon$-entropy. From the computational viewpoint, however,
the value of this quantity seems to be limited by the apparent
difficulty in establishing lower bounds. Using local approximation we will 
next define a quantity $\Rcp_L (\epsilon )$ which is
more amenable to obtaining estimates than $\Afn_L (\epsilon )$ and 
provides a ready means for obtaining upper bounds for $\Afn_L (\epsilon )$ 
(see Proposition~\ref{P-majorize}) with a view to the application
of Theorem~\ref{T-ccompact}, as we will illustrate in the case of
noncommutative tori in Example~\ref{E-tori}. 
It can also be shown (by suitably adjusting the proof of 
\cite[Prop.\ 3.9]{K} for instance) that by taking 
$\limsup_{\epsilon\to 0+}\log\Rcp_L (\epsilon ) / \log (\epsilon^{-1})$ we 
obtain a generalization of the Kolmogorov dimension of a compact metric 
space, whose utility depends on our ability to estimate 
$\log\Rcp_L (\epsilon )$ from below. We will not be concerned
here with obtaining lower bounds for $\log\Rcp_L (\epsilon )$,
but we point out that this can often be done by using the 
Hilbert space geometry implicit in the given operator system or 
$C^*$-algebra as in \cite{K}.
 
For a Lip-normed nuclear operator system $(X,L)$ and $\epsilon > 0$ we 
denote by $\CPA_L (\epsilon )$ the collection of triples
$(\alpha , \beta , B)$ where $B$ is a finite-dimensional $C^*$-algebra
and $\alpha : X\to B$ and $\beta : B\to X$ are u.c.p.\ maps with
$\| (\beta\circ\alpha )(x) - x \| < \epsilon$ for all $x\in
\mathcal{D}_1 (L)$. This collection is non-empty by 
Proposition~\ref{P-normLip-norm}. 
Conversely, if $(X,L)$ is any Lip-normed operator system $(X,L)$ and
$\CPA_L (\epsilon )$ is non-empty for each $\epsilon > 0$, then $X$ is
nuclear owing to the density of $\mathcal{D}(L)$ in $X_{\sa}$.

\begin{definition}
Let $(X,L)$ be a Lip-normed nuclear operator system. For
$\epsilon > 0$ we set
$$ \Rcp_L (\epsilon ) = \inf \{ \rank (B) : (\alpha , \beta , B)\in
\CPA_L (\epsilon ) \} . $$ 
\end{definition}

\begin{proposition}\label{P-majorize}
If $(X,L)$ is a Lip-normed nuclear operator system and $\epsilon > 0$
then $\Rcp_L (\epsilon )\geq\Afn_L (\epsilon )$.
\end{proposition}

\begin{proof}
Let $\epsilon > 0$. Then there is a
triple $(\alpha , \beta , B)\in\CPA_L (\epsilon )$ with $\rank (B) = 
\Rcp_L (\epsilon )$.
Set $Y=\alpha (B)$, and let $L_Y$ be the Lip-norm on $Y$ induced by $L$ 
via $\alpha$. Then $\dist_s (X,Y)\leq\epsilon$ by 
Proposition~\ref{P-quotientest}, and since $B$ unitally embeds into a
matrix algebra of the same rank we obtain $\Afn_L (\epsilon )\leq
\Rcp_L (\epsilon )$.
\end{proof}

\begin{example}[Noncommutative tori]\label{E-tori}
Let $\rho : \mathbb{Z}^d \times \mathbb{Z}^d \to \mathbb{T}$ be an
antisymmetric bicharacter and for $1\leq i,j \leq k$ set
$$ \rho_{ij} = \rho (e_i , e_j ) $$
with $\{ e_1 , \dots , e_d \}$ the standard basis for $\mathbb{Z}^d$.
We call the universal $C^*$-algebra $A_\rho$ generated by unitaries
$u_1 , \dots , u_d$ satisfying
$$ u_j u_i = \rho_{ij} u_i u_j $$
a {\em noncommutative $d$-torus}. 
Given a noncommutative $d$-torus $A_\rho$ with generators
$u_1 , \dots , u_d$ there is an ergodic action $\gamma : \mathbb{T}^d
\cong (\mathbb{R} / \mathbb{Z})^d \to \Aut (A_\rho )$ determined 
on the generators by
$$ \gamma_{(t_1 , \dots , t_d )}(u_j )
= e^{2\pi it_j} u_j $$
(see \cite{OPT}). Let $\ell$ be a length function on $\mathbb{T}^d$
(for instance, we could take the distance to $0$ with respect to the 
metric induced from the Euclidean metric on $\mathbb{R}^d$). 
By Example~\ref{E-actions} we then obtain a Lip-norm $L$ arising 
from the action $\gamma$ and length function $\ell$.
Let $\tau$ be the tracial state on $A_\rho$ defined by
$$ \tau (a) = \int_{\mathbb{T}^d} 
\gamma_{(t_1 , \cdots , t_d)}(a) \, d(t_1 , \dots , t_d ) $$
for all $a\in A_\rho$, where $d(t_1 , \dots , t_d )$ is normalized Haar
measure.

Let $\mathcal{A}(d,\ell )$ be the subset of $\mathcal{R}$ consisting
of all noncommutative $d$-tori Lip-normed as above with respect to the
length function $\ell$. This is in fact a subset of
$\mathcal{R}_{\fa}$ by Proposition~\ref{P-finiterank}, since
noncommutative tori are nuclear. We will show using 
Theorem~\ref{T-ccompact} that $\mathcal{A}(d,\ell )$ is totally bounded.

For $(n_1 , \dots , n_d )\in\mathbb{N}^d$ we denote by
$R(n_1 , \dots , n_d )$ the set of points $(k_1 , \dots ,
k_d )$ in $\mathbb{Z}^d$ such that $|k_i | \leq n_i$ for
$i=1, \dots , d$. For $a\in A_\rho$, we 
define for every $(n_1 , \dots , n_d )\in\mathbb{N}^d$ the partial 
Fourier sum
$$ s_{(n_1 , \dots , n_d )} (a) = \sum_{(k_1 , \dots , k_d )\in
R(n_1 , \dots , n_d )} \tau (au_d^{-k_d}\cdots
u_{1}^{-k_1}) u_{1}^{k_1}\cdots u_d^{k_d} $$
and for each $n\in\mathbb{N}$ the Ces\`{a}ro mean
$$ \sigma_n (a) = (n+1)^{-d}\sum_{(n_1 , \dots , n_d )\in R(n, n, \dots , n)}
s_{(n_1 , \dots , n_d )}(a) . $$
As in classical Fourier analysis (see for example \cite{Kat})
it can be shown that if $K_n$ is the Fej\'{e}r kernel 
$$ K_n (t) = \sum_{k=-n}^n \left( 1 - \frac{|k|}{n+1} \right) e^{2\pi ikt}
= \frac{1}{n+1} \left( \frac{\sin ((n+1)t/2 )}{\sin (t/2 )} \right)^2 . $$
then for all $a\in A_\rho$ and 
$n\in\mathbb{N}$ we have
$$ \| a - \sigma_n (a) \| \leq \sum_{k=1}^d \int_{\mathbb{T}}
\| a - \gamma_{r_k (t)}(a) \| K_n (t)\, dt $$
where $r_k (t)$ denotes the $d$-tuple which is $t$ at the $k$th
coordinate and $0$ elsewhere, and $dt$ is normalized Haar measure 
(see for example the proof of \cite[Thm.\ 22]{LADVNA}). It follows that 
if $a\in\mathcal{D}(L)$ then
$$ \| a - \sigma_n (a) \| \leq\sum_{k=1}^d \int_{\mathbb{T}}
\ell (r_k (t)) K_n (t)\, dt . $$
In particular the Ces\`{a}ro means of elements of $\mathcal{D}(L)$
converge at a rate which does not depend on $\rho$. By 
\cite[Lemma 9.4]{GHDQMS} there is a constant $M > 0$ such that 
$\diam (A_\rho ,L) \leq M$ for all $(A_\rho ,L)\in\mathcal{A}(d,\ell )$.
Hence to obtain the total boundedness of $\mathcal{A}(d,\ell )$ we need
only check condition (ii) in Theorem~\ref{T-ccompact}. Let $\epsilon > 0$.
If $B$ is a finite-dimensional $C^*$-algebra
and $\alpha : A_\rho \to B$ and $\beta : B\to A_\rho$ are u.c.p.\ maps with
$\| (\beta\circ\alpha )(x) - x \| < \epsilon$ for all $x\in
\mathcal{D}_1 (L) \cap\mathcal{B}^{A\rho}_M$, then it is readily seen that
Proposition~\ref{P-normLip-norm} implies that
the triple $(\alpha , \beta , B)$ lies in $\CPA_L (\epsilon )$.
From above there is an $n\in\mathbb{N}$ such that 
each element of $\mathcal{D}_1 (L) \cap\mathcal{B}^{A_\rho}_M$ is 
within $\epsilon$ of its $n$th Ces\`{a}ro mean, which is a linear
combination of elements in $\{ u_1^{k_1}\cdots u_d^{k_d} : 
| k_i |\leq n \}$
with coefficients bounded in modulus by $M$ (since the operation of
taking a Ces\`{a}ro mean decreases the moduli of Fourier coefficients, which
are bounded by the norm of the given element). Thus in view of
Proposition~\ref{P-majorize} it suffices to show the existence
of a finite-dimensional $C^*$-algebra $B$ and u.c.p.\ maps 
$\alpha : A_\rho \to B$ and $\beta : B\to A_\rho$ such that
$\| (\beta\circ\alpha )(x) - x \| < \epsilon$ for all
$x\in \{ u_1^{k_1}\cdots u_d^{k_d}$ : $| k_i |\leq n \}$ and the
rank of $B$ can be chosen independently of $\rho$, and this is a 
consequence of \cite[Lemma 5.1]{Voi}.

Hanfeng Li has informed me that he can show that the
map from the space of antisymmetric bicharacters on $\mathbb{Z}^d$
to $\mathcal{A}(d,\ell )$ determined by $\rho\mapsto A_\rho$ is continuous
(as Rieffel showed for quantum Gromov-Hausdorff distance in 
\cite[Thm.\ 9.2]{GHDQMS}). 
In fact, given any field of 
strongly continuous ergodic actions of a compact group on  
a continuous field of unital $C^*$-algebras over a compact metric 
space $X$, at any point of $X$ the continuity of complete distance is 
equivalent to the local constancy 
(or, equivalently, the lower semicontinuity) of the function on $X$
which records the multiplicity of the action in the fibre algebras.
This is a result of the fact that
Li (unpublished notes) has worked out a general version of Rieffel's 
result on coadjoint orbits as described in Example~\ref{E-sphere}.
\end{example}

\end{document}